\input ./style/arxiv-ba.cfg
\documentclass[ba,preprint,linksfromyear]{imsart}
\makeatletter
   \@ifpackageloaded{natbib}{}{\usepackage{natbib}}
\makeatother
\usepackage{graphicx,amstext,amssymb,bm,enumitem}

\pubyear{2015}
\volume{10}
\issue{3}
\firstpage{627}
\lastpage{664}
\doi{10.1214/14-BA893}

\startlocaldefs
\setattribute{email}{text}{}

\newcommand{\yn}{\bm{Y}_n}
\newcommand{\xn}{\bm{X}_n}
\newcommand{\en}{\bm{e}_n}

\newcommand{\bbn}{\bm{\beta}_n}
\newcommand{\bbkn}{\bm{\beta}_{k_n}}
\newcommand{\bnh}{\hat{\bm{\beta}}_n}
\newcommand{\bknh}{\hat{\bm{\beta}}_{k_n}}
\newcommand{\sig}{\sigma^2}
\newcommand{\ipn}{\bm{I}_{p_n}}
\newcommand{\gn}{\bm{\gamma}_n}
\newcommand{\gkn}{\bm{\gamma}_{k_n}}
\newcommand{\gmn}{\bm{\gamma}_{m_n}}
\newcommand{\lp}{\left(}
\newcommand{\rp}{\right)}

\newcommand{\bnhb}{\hat{\bm{\beta}}_n^{\text{B}}}

\newcommand{\bntbg}{\tilde{\bm{\beta}}_n^{\text{B}}(g)}

\newcommand{\defined}{\mathrel{\mathop:}=}
\newcommand{\bon}{\bm{\beta}_{0n}}
\newcommand{\bokn}{\bm{\beta}_{0k_n}}
\newcommand{\bomn}{\bm{\beta}_{0m_n}}
\newcommand{\sigo}{\sigma^2_0}
\newcommand{\xx}{\bm{X}_n^T\bm{X}_n}
\newcommand{\xxi}{\left(\bm{X}_n^T\bm{X}_n\right)^{-1}}
\newcommand{\xxit}{(\bm{X}_n^T\bm{X}_n)^{-1}}

\newcommand{\nxxismall}{n(\bm{X}_n^T\bm{X}_n)^{-1}}

\newcommand{\lams}{\lambda_{n,1},\ldots,\lambda_{n,p_n}}
\newcommand{\lmax}{\lambda_{\max}}
\newcommand{\lmin}{\lambda_{\min}}

\newcommand{\linf}[1]{\left|\left|#1\right|\right|_\infty}
\newcommand{\linft}[1]{||#1||_\infty}
\newcommand{\ltwo}[1]{\left|\left|#1\right|\right|_2}

\newcommand{\gnheb}{\hat{g}_n^{\text{EB}}}

\newcommand{\gknheb}{\hat{g}_{k_n}^{\text{EB}}}

\newcommand{\aspo}{\text{ a.s.$(P_0)$}}
\newcommand{\zn}{\bm{Z}_n}
\newcommand{\mun}{\bm{\mu}_n}
\newcommand{\Sigman}{\bm{\Sigma}_n}
\newcommand{\xin}{\bm{\xi}_n}
\newcommand{\zkn}{\bm{Z}_{k_n}}
\newcommand{\mukn}{\bm{\mu}_{k_n}}
\newcommand{\Sigmakn}{\bm{\Sigma}_{k_n}}
\newcommand{\xikn}{\bm{\xi}_{k_n}}
\newcommand{\tnt}{\widetilde{T}_n}
\newcommand{\given}{\;\right|\;}

\newcommand{\thon}{\theta_{0n}}
\newcommand{\thokn}{\theta_{0k_n}}
\newcommand{\thton}{\widetilde{\theta}_{0n}}
\newcommand{\lton}{\breve{\lambda}_{0n}}
\newcommand{\lti}{\breve{\lambda}_{0n}^{-1}}

\newcommand{\etaon}{\eta_{0n}}

\newtheorem*{defn}{Definition}
\newtheorem{thm}{Theorem}
\newtheorem{lem}{Lemma}
\newtheorem*{prf}{Proof}

\endlocaldefs

\begin{document}

\begin{frontmatter}
\title{Necessary and Sufficient Conditions for~High-Dimensional
Posterior Consistency under $g$-Priors}
\runtitle{Posterior Consistency under $g$-Priors}

\begin{aug}
\author[a]{\fnms{Douglas K.} \snm{Sparks}\corref{}\ead[label=e1]{dksparks@stanford.edu}},
\author[b]{\fnms{Kshitij}
\snm{Khare}\ead[label=e2]{kdkhare@stat.ufl.edu}},
\and
\author[c]{\fnms{Malay} \snm{Ghosh}\ead[label=e3]{ghoshm@stat.ufl.edu}}

\runauthor{D. K. Sparks, K. Khare, and M. Ghosh}

\address[a]{Stanford University, \printead{e1}}
\address[b]{University of Florida, \printead{e2}}
\address[c]{University of Florida, \printead{e3}}
\end{aug}

%
\begin{abstract}
We examine \emph{necessary and sufficient} conditions for posterior
consistency under $g$-priors, including extensions to hierarchical and
empirical Bayesian models. The key features of this article are that
we allow the number of regressors to grow at the same rate as the
sample size and define posterior consistency under the sup vector norm
instead of the more conventional Euclidean norm. We consider in
particular the empirical Bayesian model of \cite{GF2000}, the
hyper-$g$-prior of \cite{LPMCB2008}, and the prior considered by
\cite{ZS1980}.
\end{abstract}

%

\end{frontmatter}


\section{Introduction}

Arnold Zellner made pioneering contributions to the fields of
statistics and econometrics.
One of his works, the $g$-prior \citep{Zellner1986}, has
become a cornerstone of research in Bayesian statistics.
The $g$-prior specifies that a vector of regression coefficients is
normally distributed a~priori with some mean (typically zero) and
covariance matrix equal to a scalar multiple (typically denoted by $g$)
of the covariance matrix of the maximum likelihood estimator.
These
priors are useful for conventional hierarchical and empirical
Bayesian analysis \citep{GSJ1982} for linear
regression models, but their application extends well beyond to variable
selection 
\citep{GF2000},
Bayesian classification of high-dimensional low--sample size data \citep
{MGG2005}, and many other interesting topics
of research. The excellent article of \cite{LPMCB2008} provides
a succinct account of mixtures of $g$-priors for Bayesian variable selection.

One very important but often neglected issue in the selection of priors is
to examine the consistency of resulting posteriors
in the frequentist sense.
We will provide a formal
definition in Section~2, but in plain language, this means that as one
accumulates more and more samples, the posterior distribution of the
parameter under consideration gets closer and closer to its true value,
eventually becoming degenerate at this point in the limit.
Recently, the notion of posterior consistency has also been considered
in nonparametric settings \citep{BSW1999, GGV2000}.

In the $g$-prior model, if the number of regressors~$p$ does not vary
with~$n$, then it can easily be seen that the resulting posterior is
inconsistent if $g$ is fixed, but the problem disappears when $g\equiv
g_n$ with $g_n\to\infty.$ See Section~2 for the details of these
results. Now suppose instead that the number of regressors~$p\equiv
p_n$ increases with~$n$ but satisfies $p_n<n$ and $p_n/n\to\alpha,$
where $0\le\alpha<1.$
This situation represents the so-called ``large~$p$, large~$n$''
regime, which has been considered in the context of model selection.
\cite{BGN2003} provide scenarios where the Bayes factor is consistent
but the Bayesian Information Criterion (BIC) is not, with the
explanation that BIC may be a poor approximation to the Bayes factor
when $p_n\to\infty$. \cite{MGC2010} examine consistency of the Bayes
factor for nested normal linear models with $p_n\to\infty$, including
the case where $p_n$ grows at the same rate as the sample size.
Also, \cite{Jiang2007} addressed the variable selection problem when
$p_n>n$ and provided convergence rates for the fitted densities in a
broad class of generalized linear models.

In the context of parameter estimation as examined here, \cite
{Ghoshal1999} considered certain types of high-dimensional linear
models and provided a valuable contribution by proving not only
posterior consistency but also asymptotic normality of the posterior
distribution. However, our work differs from \cite{Ghoshal1999} in
three principal respects. First, and perhaps most fundamentally, the
$g$-prior model itself involves an unknown sampling variance~$\sig$
with an associated prior (the prior on the regression coefficients is
taken to be conditional on $\sig$). Such a structure is not included in
the class of models considered in \cite{Ghoshal1999}. Second, we
provide \emph{necessary and sufficient} conditions for posterior
consistency in three of the four $g$-prior models we consider. While we
readily admit that stronger results such as asymptotic normality are
perhaps more useful whenever posterior consistency occurs, our \emph
{necessary} conditions demonstrate circumstances in which posterior
consistency fails to occur at all, which we believe to be interesting
in their own right. Third, our work allows the parameter space for the
$p_n$-dimensional vector of regression coefficients to be taken as
$\mathbb{R}^{p_n}$, as is natural. This contrasts with \cite
{Ghoshal1999}, which essentially requires the restriction of the
parameter space to a sequence of compact sets.

\cite{Bontemps2011} also extended the work of \cite{Ghoshal1999} in
several ways by permitting the model to be misspecified and the number
of regressors to grow proportionally to the sample size, the latter of
which is also a feature of our work. However, our work differs from
\cite{Bontemps2011}, most notably by allowing the consideration of
models where the sampling variance $\sigma^2$ is assumed to be unknown.
There are also differences in the assumptions. In particular, \cite
{Bontemps2011} does not make any assumption analogous to the eigenvalue
bounds that we will later impose in~(A\ref
{item:eigenvalue-assumption}). On
the other hand, unlike \cite{Bontemps2011}, we do not make any
assumptions on the asymptotic behavior of
the true coefficient vector~$\bon$. We must also emphasize once again
that, unlike \cite{Bontemps2011}, we provide conditions that are both
\emph{necessary and sufficient} for posterior consistency. This
establishes circumstances in which posterior consistency definitively
does \emph{not} occur, which can in some cases be rather surprising
(see the remarks following Theorem~\ref{thm:pcg-eb}, for example). The
recent work of
\cite{ADLBS2013} establishes sufficient conditions for posterior
consistency in linear models under
shrinkage priors. Again, the most notable difference between the models
considered in
\cite{ADLBS2013} and the $g$-prior based models considered in this
paper is that the variance parameter $\sigma^2$ is assumed to be known
in \cite{ADLBS2013}.
\cite{LO2013} consider a high dimensional Bayesian Principal Components
Analysis regression setup with $p_n > n$ and normal priors, and examine
posterior consistency (in the $\ell_2$-norm) and convergence rates under
appropriate assumptions on the rank of the design matrix.

Other authors have addressed the asymptotic properties of $g$-prior
models, but for model selection instead of parameter estimation. \cite
{FLS2001} provided both
theoretical results and simulation-based evidence for the consistency
of posterior model probabilities under particular choices for the
$g$-prior hyperparameter $g\equiv g_n$. \cite{LPMCB2008} took a more
theoretical approach and proved the consistency of posterior model
probabilities under hierarchical and empirical Bayesian $g$-prior
models, but only in the case where the dimensionality~$p_n$ of the full
model is fixed. More recently, \cite{SC2011} provided similar results
in the case where $p_n\to\infty$, albeit under a considerable number of
assumptions. They also note that these results can be extended to
$p_n>n$, the so-called ``large~$p$, small~$n$'' regime, when combined
with certain dimension reduction approaches. See also the work of \cite
{ZJY2009}.

Another new feature of our work is that we have established posterior
consistency under the sup vector
norm $\ell_\infty\; (||\bm{x}||_\infty= \max_{1\le i\le p}|x_i|)$
rather than the conventional $\ell_2 \;
(||\bm{x}||_2 = [\sum_{i=1}^p x_i^2]^{1/2})$ vector norm. The choice is
motivated primarily because the $\ell_\infty$~norm introduces added
flexibility to our procedure, since it is weaker than the $\ell_2$~norm
(as a vector norm), noting that $||\bm{x}||_\infty\le||\bm{x}||_2.$ In
particular, for proving consistency when the number of covariates $p_n$
grows with the sample size, the sup norm approach allows $p_n$ to grow
at a faster rate than is possible under the $\ell_2$ norm. The simplest
yet most convincing fact in this
regard is the following. For the linear model $\bm{Y}_n = \bm{X}_n \bbn
+ \bm{e}_n$ with i.i.d.\ Gaussian errors and $\bm{X}_n^T \bm{X}_n = n
{\bf I}_{p_n}$ (orthogonal covariates, ${\bf I}_{p_n}$ denotes the
identity matrix of dimension $p_n$),
the MLE for $\bbn$ is consistent under the $\ell_2$ vector norm if and
only if $p_n = o(n)$. However, the MLE for $\bbn$ is
still consistent under the $\ell_\infty$~norm for any $p_n < n$. See
remark immediately following Lemma 1.

As discussed above, if $p_n \rightarrow\infty$, it is harder to prove
posterior consistency under the
$\ell_2$~norm as compared to the $\ell_\infty$~norm. However, in the
same vein, it is harder to prove
posterior inconsistency under the $\ell_\infty$~norm as compared to the
$\ell_2$~norm. In particular,
any necessary condition for posterior consistency under the $\ell_\infty
$~norm is also a necessary
condition for posterior consistency under the $\ell_2$~norm. Hence,
this paper also provides novel
necessary conditions for posterior consistency under the conventional
$\ell_2$~norm (note that
assumption (A2) in Section $2$ subsumes the case $p_n/n \rightarrow0$).

The outline of the remaining sections is as follows. Section~2 provides
necessary and sufficient conditions for posterior consistency for a
nonstochastic sequence $\{g_n,n\ge1\}.$ In the process, we demonstrate
the posterior consistency or inconsistency of some popular
recommendations regarding the choice of $g_n.$ Section~3 provides
necessary and sufficient conditions for posterior consistency in an
empirical Bayesian context in which $g_n$ is estimated from the data.
Section~4 provides necessary and sufficient conditions for posterior
consistency under the hierarchical hyper-$g$-prior model \citep
{LPMCB2008}. Section~5 considers the celebrated Zellner-Siow prior
\citep{ZS1980} and provides a sufficient (though not necessary)
condition for posterior consistency under this model. At the end of
each of Sections 2--5, the interpretations and implications of the
results are briefly discussed. Some final remarks are made in Section~6.
It should be noted that although the key results of Sections~3--5 yield
the same condition for posterior consistency, the techniques used to
prove these results differ substantially among the three models.
Furthermore, the coincidence of the conditions in Theorems~2--4 should
not be misconstrued as a suggestion that the same condition would be
shared by other hierarchical or empirical Bayesian $g$-prior models.
Specifically, this condition is not shared by Theorem~1, yet the
non-hierarchical model addressed by Theorem~1 can be considered as a
hierarchical model with a sequence of degenerate hyperpriors. Moreover,
it should again be noted that the conditions in Theorems~1--3 are both
necessary and sufficient, but the condition for the Zellner-Siow
$g$-prior model provided in Theorem~4 is merely sufficient, and its
necessity or lack thereof is not presently clear.


\section{Non-Hierarchical Model}

Consider the usual linear model $\yn=\xn\bbn+\en,$
with response $\yn=(Y_{n,1},\ldots,Y_{n,n})^T$,
covariates $\xn=(\bm{x}_{n,1},\ldots,\bm{x}_{n,n})^T$, regression
coefficients $\bbn=(\beta_{n,1},\ldots,\beta_{n,p_n})^T$ and errors $\en
=(e_1,\ldots,e_n)^T$.
We now impose the following assumptions:
\renewcommand{\labelenumi}{(A\arabic{enumi})}
\begin{enumerate}[leftmargin=*]
\item\label{item:normality-assumption} The errors are distributed as
$\en\sim N_n(\bm{0}_n,\sig{\bf I}_n)$. Here
$\bm{0}_n$ denotes the vector of length $n$ with all zero entries.
\item The number of regressors~$p_n$ is a nondecreasing sequence with
$p_n<n$ and $p_n/n\to\alpha$, where $0\le\alpha<1.$
\item\label{item:eigenvalue-assumption} The eigenvalues $\lams$ of the
matrix $\nxxismall$ satisfy
$0<\lmin\le\inf_{n,i}\lambda_{n,i}\le\sup_{n,i}\lambda_{n,i}\le\lmax
<\infty$ for some $\lmin$ and $\lmax$.
\end{enumerate}
Note that (A\ref{item:eigenvalue-assumption}) implies that $\lmax
^{-1}\ipn\le n^{-1}\bm{X}_n^T\bm{X}_n\le\lmin^{-1}\ipn.$ This
assumption is identical to assumption (A2) of
\cite{ADLBS2013}.

The goal in such a model is estimation of $\bbn.$ Minimal sufficiency
leads to the reduction $(\bnh,S_n),$ where $\bnh=\xxit\xn^T\yn,$ the
maximum likelihood estimator of $\bbn,$ and $S_n=||\yn-\xn\bnh||_2^2,$
the error sum of squares.
Note that conditional on $\bbn$ and $\sig,$ $\bnh$ and $S_n$ are
mutually independent with $\bnh\mid\bbn,\sig\sim N_{p_n}(\bbn,\sig\xxit
)$ and $S_n\mid\bbn,\sig\sim\sig\chi^2_{n-p_n}.$

Now suppose priors are specified as $\bbn\mid\sig\sim N_{p_n}(\gn,g\sig
\xxit)$ (Zellner's $g$-prior) and $\sig\sim\text
{InverseGamma}(a/2,b/2),$ where we permit $a\ge-2$ and $b\ge0$ to
accommodate such improper priors as $\pi(\sig)\propto1/\sig,1/\sigma
,\text{ or }1.$
Suppose further that $g\equiv g_n$ is specified as a known sequence of
constants.
This collection of likelihoods and priors comprises our
non-hierarchical $g$-prior model, which we denote by $P_M.$ One
motivation for the use of such a model is the convenient form of the
Bayes estimator under squared error loss,
\[
\bnhb\defined E_M(\bbn\mid\bnh,S_n)=\frac{g_n}{g_n+1}\bnh+\frac
{1}{g_n+1}\gn,
\]
where $\bnh$ denotes the MLE.

We now introduce the formal definition of posterior consistency.

\begin{defn}
Let $\bon\in\mathbb{R}^{p_n}$ for each $n\ge1$, and let $\sigo>0.$ Now
let $P_0$ denote the distribution of $\{(\bnh,S_n),n\ge1\}$ under the
model $\yn=\xn\bon+\en$, where $\en\sim N_n(\bm{0}_n,\sigo I_{n})$, for
each $n\ge1.$ The sequence of posterior distributions $P_M(\bbn\mid\bnh
,S_n)$ is said to be consistent under the $\ell_\infty$~norm at $\{(\bon
,\sigo),n\ge1\}$ if $P_M(||\bbn-\bon||_\infty>\epsilon\mid\bnh,S_n)\to
0\,\aspo$ for every $\epsilon>0.$
\end{defn}

It should be immediately noted that the type of posterior consistency
considered herein is fundamentally different from what could instead be
considered in the analysis of Bayesian methodology, that is,
convergence of the posterior under the same model~$P_M$ under which it
is derived. In this case, one is assuming that the prior associated
with the model~$P_M$ is in some sense ``true.'' However, this approach
is perhaps too favorable in that posterior consistency is quite easy to
achieve. In fact, in this approach, a quite general result due to \cite
{Doob1948} states that posterior consistency occurs on a set of
parameter values with probability~1 under the prior associated with
$P_M.$ Instead, the type of posterior consistency considered herein is
fundamentally frequentist in nature, that is, the values~$\bon$ and
$\sigo$ are considered fixed but unknown.

The frequentist properties of Bayesian methods have been of interest
for some time. Even pure frequentists may be interested in originally
Bayesian procedures, or limits and approximations thereof, due to
considerations such as admissibility and the convenient elimination of
nuisance parameters. Indeed, it was shown as early as \cite{Laplac1774}
that in simple cases, the posterior distribution and the distribution
of the maximum likelihood estimator are comparable for large sample
sizes. More sophisticated versions of such results have been developed
in more recent times \citep{Bernstein1934, DF1986, GSJ1982, LeCam1982,
vonMises1964}.

We now provide a lemma establishing strong frequentist consistency of
the MLE $\bnh$ in the $\ell_\infty$ norm.

%
\begin{lem}
\label{lem:mle-linf}
Let $\bm{Z}_n\sim N_{p_n}(\bm{0}_{p_n},n^{-1}\bm{V}_n),$ where $p_n<n$,
and where the eigenvalues $\omega_{n,1},\ldots,\omega_{n,p_n}$ of $\bm
{V}_n$ satisfy $\sup_{n,i}\omega_{n,i}=\omega_{\max}<\infty.$
Then $||\bm{Z}_n||_\infty\to0\,\text{ almost surely}.$
\end{lem}
%
%
\begin{prf}[Proof of Lemma~\ref{lem:mle-linf}]
First note that $\text{Var}(Z_{n,i})=n^{-1}V_{n,ii}\le n^{-1}\omega
_{\max}$, and $n^{1/2}V_{n,ii}^{-1/2}Z_{n,i}\sim N(0,1).$ Now let
$\epsilon>0$. Since
\begin{align*}
\sum_{n=1}^\infty P\lp||\bm{Z}_n||_\infty>\epsilon\rp
&=
\sum_{n=1}^\infty P\lp\max_{1\le i\le p_n}\left|Z_{n,i}\right|>\epsilon
\rp,
\end{align*}
it follows that
\begin{align*}
\sum_{n=1}^\infty P\lp||\bm{Z}_n||_\infty>\epsilon\rp
&\le
\sum_{n=1}^\infty\sum_{i=1}^{p_n} P\lp n^{1/2}V_{n,ii}^{-1/2}\left
|Z_{n,i}\right|>\epsilon\lp n^{-1}V_{n,ii}\rp^{-1/2}\rp\\
&\le
\sum_{n=1}^\infty\sum_{i=1}^{p_n} P\lp n^{1/2}V_{n,ii}^{-1/2}\left
|Z_{n,i}\right|>\epsilon\, \omega_{\max}^{-1/2}n^{1/2}\rp\\
&\le
\sum_{n=1}^\infty\sum_{i=1}^{p_n} \frac{15\omega_{\max}^3}{\epsilon^6
n^3} <\infty
\end{align*}
by applying Markov's inequality to $n^{3}V_{n,ii}^{-3}Z_{n,i}^6$. The
result follows from the Borel-Cantelli lemma, noting that $p_n<n.$
\end{prf}
%

Observe that Lemma~\ref{lem:mle-linf} under $P_0$ with assumptions
(A\ref{item:normality-assumption})--(A\ref{item:eigenvalue-assumption})
and $\bm{Z}_n=\bnh-\bon$ implies that $||\bnh-\bon||_\infty\to0\,\aspo
.$ Thus, the MLE~$\bnh$ retains strong frequentist consistency in the
$\ell_\infty$ norm even as $p_n$ grows at a rate exactly proportional
to $n.$ To contrast this with the behavior of the MLE under the
conventional $\ell_2$ vector norm, note that we have the upper bound
\[
||\bnh-\bon||_2^2\le\frac{\sigo\lmax}{n}\lp\bnh-\bon\rp^T\frac{1}{\sigo
}\bm{X}_n^T\bm{X}_n\lp\bnh-\bon\rp
\]
and a similar lower bound with $\lmax$ replaced by $\lmin.$ Since
\[
\lp\bnh-\bon\rp^T\frac{1}{\sigo}\bm{X}_n^T\bm{X}_n\lp\bnh-\bon\rp\sim
\chi^2_{p_n},
\]
it can be immediately seen that $p_n=o(n)$ is required for strong
frequentist consistency of the MLE~$\bnh$ under the $\ell_2$~norm.

In a Bayesian analysis,
Lemma~\ref{lem:mle-linf} leads to the following useful lemma, which
essentially states that $\bon$ may be replaced by $\bnh$ in the
definition of posterior consistency.

%
\begin{lem}
\label{lem:pcg-basic}
In the $g$-prior model (both hierarchical and non-hierarchical),
$P_M(||\bbn-\bon||_\infty>\epsilon\mid\bnh,S_n)\to0\,\aspo$
for every $\epsilon>0$ if and only if
$P_M(||\bbn-\bnh||_\infty>\epsilon\mid\bnh,S_n)\to0\,\aspo$ for every
$\epsilon>0.$
\end{lem}
%
%
\begin{prf}[Proof of Lemma~\ref{lem:pcg-basic}]
The triangle inequality implies that
\begin{align*}
&P_M\lp\linf{\bbn-\bnh}>2\epsilon\mid\bnh,S_n\rp-
P_M\lp\linf{\bnh-\bon}>\epsilon\mid\bnh,S_n\rp\notag\\
&\qquad\le
P_M\lp\linf{\bbn-\bon}>\epsilon\mid\bnh,S_n\rp\notag\\
&\qquad\le
P_M\lp\linf{\bbn-\bnh}>\epsilon/2\mid\bnh,S_n\rp
+P_M\lp\linf{\bnh-\bon}>\epsilon/2\mid\bnh,S_n\rp.
\end{align*}
When conditioning on $\bnh$ and $S_n,$
\[
P_M\lp\linft{\bnh-\bon}>\epsilon\mid\bnh,S_n\rp=I\lp\linft{\bnh-\bon
}>\epsilon\rp,
\]
where $I(\cdot)$ denotes the usual indicator function.
Lemma~\ref{lem:mle-linf} implies that $\linft{\bnh-\bon}\to0\,\aspo,$
from which it follows that $I(\linft{\bnh-\bon}>\epsilon)\to0\,\aspo$
for all $\epsilon>0.$ This and the above inequalities immediately yield
the result.
\end{prf}
%

To establish results on posterior consistency or inconsistency in the
non-hierarchical $g$-prior model, we first define $T_n\defined(\bnh-\gn
)^T\xx(\bnh-\gn),$ so that $T_n/\sig$ is the usual frequentist
likelihood ratio test statistic for a test of $H_0:\bbn=\gn$ vs.\
$H_a:\bbn\ne\gn$ with known variance~$\sig.$
Then the joint posterior $\pi_n(\bbn,\sig\mid\bnh,S_n)$ is given by
\begin{align*}
\pi_n(\bbn,\sig\mid\bnh,S_n)&\propto
\exp\left[-\frac{1}{2}\lp\bbn-\bnhb\rp^T\lp\frac{g_n}{g_n+1}\sig\xxi\rp^{-1}
\lp\bbn-\bnhb\rp\right]\notag\\
&\qquad\qquad\times
\lp\sig\rp^{-(n+p_n+a)/2}\exp\left[-\frac{1}{2\sig}\lp S_n+b+\frac
{T_n}{g_n+1}\rp\right],
\end{align*}
and integrating out $\bbn$ from this yields the marginal posterior of
$\sig,$
\begin{align*}
\pi_n(\sig\mid\bnh,S_n)
&\propto
\lp\sig\rp^{-(n+a)/2}\exp\left[-\frac{1}{2\sig}\lp S_n+b+\frac
{T_n}{g_n+1}\rp\right],
\label{eq:post-s}
\end{align*}
i.e., $\sig\mid\bnh,S_n\sim\text{InverseGamma}((n+a-2)/2,\;\tnt/2),$
where we define $\tnt\defined S_n+b+(g_n+1)^{-1}T_n.$
For notational convenience, for each~$n\ge1,$ define
\begin{align*}
\lton&\defined
\frac{n||\gn-\bon||_2^2}{(\gn-\bon)^T\xx(\gn-\bon)},\\
\thon&\defined E_0(T_n)=p_n\sigo+n\lti\ltwo{\gn-\bon}^2, \\
\thton&\defined E_0(\tnt)=(n-p_n)\sigo+b+\frac{1}{g_n+1}\lp p_n\sigo
+n\lti\ltwo{\gn-\bon}^2\rp,
\end{align*}
and note that $\lmin\le\lton\le\lmax$ since $\lmax^{-1}\ipn\le n^{-1}\xx
\le\lmin^{-1}\ipn.$

The following lemmas establish the behavior of various quantities under
$P_0,$ and they will be heavily used in proving posterior consistency
or inconsistency in both the non-hierarchical and hierarchical
$g$-prior models. The proof of each lemma can be found in the Appendix.

%
\begin{lem}
\label{lem:sn}
$(n-p_n)^{-1}S_n\to\sigo\,\aspo.$
\end{lem}
%

%
\begin{lem}
\label{lem:tn-to-mean}
If $\alpha>0$ or $\liminf_{n\to\infty}||\gn-\bon||_2^2>0,$ then
$T_n/\thon\to1\,\aspo.$
\end{lem}
%

%
\begin{lem}
\label{lem:tnt-to-mean}
$\tnt/\,\thton\to1\,\aspo.$
\end{lem}
%

The following lemmas regarding the normal distribution will be useful
in establishing the condition for posterior consistency in the
non-hierarchical case. The proofs are provided in the Appendix.

%
\begin{lem}
\label{lem:normal-means-noncvgc}
Let $\zn\sim N_{p_n}(\mun,\Sigman),\;\Sigman$ positive definite, $n\ge1.$
If $\linf{\mun-\xin}\nrightarrow0,$ then there exist $\epsilon>0$ and a
subsequence $k_n$ of $n$ such that
$P(\linf{\zkn-\xikn}>\epsilon)\ge1/2$ for all~$n.$
\end{lem}
%

%
\begin{lem}
\label{lem:normal-center}
Let $Z\sim N(\mu,\tau^2).$ Then $P(|Z|\le\xi)\le1-2\Phi(-\xi/\tau)$
for every $\xi\ge0,$
where $\Phi$ is the standard normal cdf.
\end{lem}
%

%
\begin{lem}
\label{lem:normal-correlations}
Let $\zn\sim N_{p_n}(\bm{0}_{p_n},\bm{\Sigma}_n)$ for each $n\ge1,$
where $\bm{\Sigma}_n$ has each diagonal entry equal to~1 and
eigenvalues $\omega_{n,1},\ldots,\omega_{n,p_n}.$ If $\inf_{n,i}\omega
_{n,i}=\omega_{\min},$ then
$\inf_{n,i}\text{Var}(Z_i\mid Z_{i+1},\ldots,Z_{p_n})\ge\omega_{\min}.$
\end{lem}
%

Finally, one additional lemma provides a key result about the marginal
posterior of~$\sig.$ Again, the proof is deferred to the Appendix.

%
\begin{lem}
\label{lem:sig-bounds}
In the non-hierarchical $g$-prior model, the posterior distribution
of~$\sig$ satisfies $P_M(\thton/2n\le\sig\le2\thton/n\mid\bnh,S_n)\to
1\aspo.$
\end{lem}

Note that although $g_n$ does not appear explicitly in the result in
Lemma~\ref{lem:sig-bounds}, the result nevertheless does depend on the
choice of $g_n$ since it is involved in the quantity $\thton.$

We now state and prove the necessary and sufficient condition for
posterior consistency in the non-hierarchical $g$-prior model.

%
\begin{thm}
\label{thm:fixed-g}
In the non-hierarchical $g$-prior model $P_M,$ posterior consistency
occurs if and only if both $(g_n+1)^{-1}||\gn-\bon||_\infty\to0$ and
$g_n(g_n+1)^{-2}(\log p_n)n^{-1}||\gn-\bon||_2^2\to0.$
\end{thm}

\noindent
The proof of this theorem in provided in the Appendix.

\subsection{Interpretations and Implications}
In the same vein as frequentist consistency, posterior consistency can
be conceptualized as the idea that the center (not necessarily the
mean) of the posterior distribution converges to the true value while
the spread (not necessarily the variance) of the posterior distribution
converges to zero. In light of this, it is noteworthy that the two
conditions in Theorem~\ref{thm:fixed-g} arise from precisely such
considerations. The first condition controls the convergence to zero of
the $\ell_\infty$-distance between the posterior's center and the true
value $\bon,$ while the second condition controls the convergence of
the posterior's spread to zero. Both conditions are necessary for
posterior consistency to hold.

In the simple case where $p_n$ does not increase with $n,$ it is
typical to fix the prior mean as $\gn=\bm{\gamma}$ and to assume that
$\bm{\beta}_{0n}=\bm{\beta}_0$ also does not vary with~$n$. In this
case it can be immediately seen that although the second condition of
Theorem~\ref{thm:fixed-g} is satisfied, the first condition fails
except in the serendipitous case that $\bm{\gamma}=\bm{\beta}_0.$ Of
course, the result is somewhat obvious even without appealing to
Theorem~\ref{thm:fixed-g}, since the posterior mean is simply a
weighted average of the MLE~$\hat{\bm{\beta}}_n,$ which is strongly
consistent for $\bm{\beta}_0,$ and the prior mean~$\bm{\gamma}$ with
weights~$g(g+1)^{-1}$ and~$(g+1)^{-1}.$ In this case, the situation may
be remedied by taking any choice of $g_n$ that tends to infinity. For
instance, the unit information prior \citep{KW1995} is equivalent to
taking $g_n=n,$ while $g_n=\max\{n,p_n^2\}$ has also been recommended
\citep{FLS2001}. Either choice yields posterior consistency in the
fixed-$p$ case.

The result of Theorem~\ref{thm:fixed-g} becomes more interesting when
$p_n\to\infty.$ 
Suppose that $||\gn-\bon||_\infty=O(1),$ but $||\gn-\bon||_2^2=O(p_n)$.
This can happen, for example, if (a) ${\boldsymbol\gamma}_n = \hat{\bm
{\beta}}_n$, or (b) the entries
of $\bon$ are uniformly bounded and $\bm{\gamma}_n = c \hat{\bm{\beta
}}_n$ where $0 \leq c < 1$ (follows immediately from Lemma \ref{lem:mle-linf}).
In this case, the first condition is satisfied as long as $g_n\to\infty
,$ but the second condition imposes the additional requirement that
$g_n$ must grow faster than $p_n n^{-1}\log p_n.$ The aforementioned
choices of $g_n=n$ or $g_n=\max\{n,p_n^2\}$ provide posterior
consistency in this case as well.

As another special case, suppose $p_n=O(n)$ exactly, but suppose only a
finite number $m>0$ of components of $\gn-\bon$ are nonzero and these
$m$~components remain fixed as $n$ grows. This circumstance could arise
with the logical choice $\gn=\bm{0}_{p_n}$ if only the first few
covariates are present in the ``true'' frequentist model~$P_0,$ but
covariates continue to be added as the sample size increases. Then any
$g_n\to\infty$ ensures posterior consistency. This case is admittedly
uninteresting in the non-hierarchical model, but we will revisit its
behavior later under empirical and hierarchical Bayesian models.



\section{Empirical Bayesian Model}

A popular approach is to avoid specifying $g$ or $g_n$ altogether by
the use of an empirical Bayes method (George and Foster, 2000) in which
the value of $g$ is estimated from the data. The most common technique
is to use the value of $g$ that maximizes its marginal likelihood,
restricted to $g\ge0.$ By integrating out $\bbn$ and $\sig$ from the
joint distribution of $\bnh,S_n,\bbn,\sig,$ the marginal likelihood of
$g$ is found to be
\[
L(g;\bnh,S_n)\propto(g+1)^{(n-p_n+a-2)/2}\big[(g+1)(S_n+b)+T_n\big
]^{-(n+a-2)/2},
\]
for which the maximizing value of $g$ subject to $g\ge0$ is
\[
\hat{g}_n^{\text{EB}}\defined\max\left\{0,\lp\frac{n-p_n+a-2}{S_n+b}\rp
\lp\frac{T_n}{p_n}\rp-1\right\}.
\]
We first provide a lemma (proven in the Appendix) that addresses the
behavior of $\gnheb.$

%
\begin{lem}
\label{lem:eb-g-bound}
If $\liminf_{n\to\infty}||\gn-\bon||_2^2>0,$ then
$\liminf_{n\to\infty}\gnheb>0\,\aspo.$
\end{lem}
%


Since $\hat{g}_n^{\text{EB}}$ is simply a function of $(\bnh,S_n),$ the
empirical Bayes posterior is identical to the simple non-hierarchical
Bayes posterior, but with the data-dependent quantity $\hat{g}_n^{\text
{EB}}$ in place of $g_n.$ Thus, while Theorem~\ref{thm:fixed-g} would
allow us to immediately state a necessary and sufficient condition for
posterior consistency in terms of $\gnheb,$ an alternative condition
not involving data-dependent quantities would be preferable. The
following result gives precisely such a condition and establishes its
necessity and sufficiency.

%
\begin{thm}
\label{thm:pcg-eb}
In the empirical Bayes $g$-prior model,
posterior consistency occurs if and only if either $\alpha=0$ or there
does not exist a subsequence~$k_n$ of $n$ and a constant $A>0$ such
that $||\gkn-\bokn||_2^2\to A$ and $||\gkn-\bokn||_\infty\nrightarrow0.$
\end{thm}

The proof of this theorem is provided in the Appendix.

\subsection{Interpretations and Implications}

It should be noted that there is no immediately obvious remedy for
inconsistency in an empirical Bayesian $g$-prior model due to the
failure of the conditions in
Theorem~\ref{thm:pcg-eb}. For any particular non-hierarchical $g$-prior
model, Theorem~\ref{thm:fixed-g} implies that there always exists a
choice of $g_n$ growing
sufficiently fast to ensure posterior consistency (although the choice
may depend on $\bon$). However, such options are not available in the
empirical Bayes
approach, since $g$ is selected via a specified function of the data.

Another salient consequence of Theorem~\ref{thm:pcg-eb} is that if
$p_n=o(n),$ then the empirical Bayes model exhibits posterior
consistency for all values of $\gn$ and $\bon.$
However, if $p_n=O(n)$ exactly, then the situation is not as simple.
For example, if $\bm{\gamma}_n = \bm{0}_{p_n}$ for every $n$ and $\lim
_{n \rightarrow\infty}
\|\bm{\beta}_{0n}\|_2^2 = \infty$, then $\|\bm{\gamma}_{k_n} - \bm{\beta
}_{0k_n}\|_2^2$ converges to $\infty$ for every subsequence $k_n$,
which implies that posterior consistency
occurs. Similarly, if $\bm{\gamma}_n = \hat{\bm{\beta}}_n$ for every
$n$, then by Lemma \ref{lem:mle-linf}, $\|\bm{\gamma}_{k_n} - \bm{\beta
}_{0k_n}\|_\infty$ converges to zero
for every subsequence $k_n$, which implies that posterior consistency
occurs. On the other hand, suppose that only a fixed number $p^\star>0$
of components of $\gn-\bon$ are
nonzero and these $p^\star$~components remain fixed as $n$ grows. Then
clearly both $||\gn-\bon||_\infty$ and $||\gn-\bon||_2^2$ converge to
constants, so the condition of
Theorem~\ref{thm:pcg-eb} fails, and the posterior is inconsistent.

This behavior is perhaps somewhat surprising. If the prior mean~$\gn$
is imagined as a guess for the true~$\bon,$ then one might speculate
that posterior inconsistency would only
occur when the guess is quite bad, i.e., when $||\gn-\bon||_2^2$ or
$||\gn-\bon||_\infty$ grows too quickly. However, in the empirical
Bayesian setting, Theorem~\ref{thm:pcg-eb} shows
that this is not the case. Intuitively, the reason is that if we allow
the data to determine the value of $g,$ then a prior mean~$\gn$ that is
``too close'' to $\bon$ (in the $\ell_2$ sense)
may cause the data to choose $g$ values that tend to a finite constant,
rather than to infinity, which leads to posterior inconsistency.
An open question regarding this behavior is whether this interesting
behavior is in some way dependent on the Gaussian tails imposed by the
$g$-prior model. However, the derivation
of a similar condition for a hierarchical $g$-prior model considered
later in Theorem~\ref{thm:pcg-conj} casts doubt on this possibility,
since the hierarchical model simply corresponds to
some marginal prior with heavier tails.


\section{Hyper-$g$-Prior Hierarchical Model}

An alternative approach to the specification of~$g$ is a hierarchical
model in which $g$ is considered a hyperparameter and is given a
hyperprior $\pi_n(g).$ Under this model, the joint posterior 
is given by
\begin{align*}
&\qquad\pi_n(\bbn,\sig,g\mid\bnh,S_n)\\
&\propto\exp\left[-\frac{1}{2}\left[\bbn-\bntbg\right]^T\lp\frac
{g}{g+1}\sig\xxi\rp^{-1}
\left[\bbn-\bntbg\right]\right] \notag\\
&\qquad\qquad
\times\lp\sig\rp^{-(n+p_n+a)/2}\exp\left[-\frac{1}{2\sig}\lp
S_n+b+\frac{T_n}{g+1}\rp\right]g^{-p_n/2}\,\pi_n(g),
\end{align*}
where $\bntbg\defined E(\bbn\mid g,\sig,\bnh,S_n)=g(g+1)^{-1}\bnh
+(g+1)^{-1}\gn.$
Integrating out $\bbn$ and subsequently $\sig$ yields the marginal posteriors
\begin{align}
\label{eq:post-sg}
\pi_n(\sig,g\mid\bnh,S_n)
&\propto
\lp\sig\rp^{-(n+a)/2}\exp\left[-\frac{1}{2\sig}\lp S_n+b+\frac
{T_n}{g+1}\rp\right]
(g+1)^{-p_n/2}\,\pi_n(g),\\
\pi_n(g\mid\bnh,S_n)
&\propto
(g+1)^{-p_n/2}\lp S_n+b+\frac{T_n}{g+1}\rp^{-(n+a-2)/2}\pi_n(g). \label
{eq:post-g}
\end{align}
The following technical lemma, which is proven in the Appendix,
establishes a relationship between posterior consistency in the
hierarchical $g$-prior model and the convergence of a particular
sequence of posterior probabilities. Note that the lemma makes no
assumptions on the particular form of the hyperprior~$\pi_n(g).$

%
\begin{lem}
\label{lem:pcg-hier}
In a hierarchical $g$-prior model, suppose that $n^{-3}\,T_n^2\,
E_M[g^2(g+1)^{-4}\mid\bnh,S_n]\to0\,\aspo.$ Then posterior consistency
occurs if and only if $P_M[(g+1)^{-1}||\gn-\bon||_\infty>\epsilon\mid
\bnh,S_n]\to0\,\aspo$ for every $\epsilon>0.$
\end{lem}
%

The form of the marginal posterior of $g$ in (\ref{eq:post-g}) suggests
that a convenient choice of hyperprior is $\pi_n(g)\propto(g+1)^{-c/2}$
for some constant~$c$, called the hyper-$g$-prior \citep{LPMCB2008}.
This prior is proper for $c>2,$ and there exists an argument \citep{LPMCB2008}
for taking $2<c\le4$, but we instead permit~$c$ to take any
real value in the present analysis. The hyper-$g$-prior yields the posterior
\begin{align}
\pi_n(g\mid\bnh,S_n)
&\propto
(g+1)^{-(p_n+c)/2}\lp S_n+b+\frac{T_n}{g+1}\rp^{-(n+a-2)/2} \notag\\
&\propto
(g+1)^{(n-p_n+a-c-2)/2}\big[(g+1)(S_n+b)+T_n\big]^{-(n+a-2)/2}.\label
{eq:post-g-conj}
\end{align}
It will also be useful to define the transformation
\begin{align}
u\defined\frac{(g+1)(S_n+b)}{(g+1)(S_n+b)+T_n},\qquad W_n\defined\frac
{S_n+b}{S_n+b+T_n}, \label{eq:transform}
\end{align}
so that $g\ge0$ if and only if $u\ge W_n.$
The next lemma asserts that Lemma~\ref{lem:pcg-hier} applies with this
choice of hyperprior. The proof can be found in the Appendix.

%
\begin{lem}
\label{lem:pcg-hier-conj}
With the hyper-$g$-prior,
$n^{-3}\,T_n^2\,E_M[g^2(g+1)^{-4}\mid\bnh,S_n]\to0\,\aspo.$
\end{lem}
%

To examine the behavior of the posterior probabilities in Lemma~\ref
{lem:pcg-hier} under the hyper-$g$-prior, we begin by using the
posterior in~(\ref{eq:post-g-conj}) to write
\begin{align*}
&P_M\left[\left.\lp\frac{1}{g+1}\rp\linf{\gn-\bon}>\epsilon\;\right|\;
\bnh,S_n\right] \\
&\qquad=P_M\left[\left.g<\frac{1}{\epsilon}\linf{\gn-\bon}-1\;\right|\;
\bnh,S_n\right] \\
&\qquad=
\frac{
\displaystyle\int_0^{q_n(\epsilon)}
(g+1)^{(n-p_n+a-c-2)/2}\big[(g+1)(S_n+b)+T_n\big]^{-(n+a-2)/2}\;dg
}
{
\displaystyle\int_0^\infty
(g+1)^{(n-p_n+a-c-2)/2}\big[(g+1)(S_n+b)+T_n\big]^{-(n+a-2)/2}\;dg
},
\end{align*}
where we define $q_n(\epsilon)\defined\max\{0,\,\epsilon^{-1}||\gn-\bon
||_\infty-1\}.$
Now define
\[
\widetilde{L}_n(\epsilon)\defined\frac{\epsilon^{-1}\linf{\gn-\bon
}(S_n+b)}{\epsilon^{-1}\linf{\gn-\bon}(S_n+b)+T_n},
\qquad L_n(\epsilon)\defined\max\left\{W_n,\widetilde{L}_n(\epsilon
)\right\},
\]
and apply the transformation in (\ref{eq:transform}) to obtain
\begin{align}
&P_M\left[\left.\lp\frac{1}{g+1}\rp\linf{\gn-\bon}>\epsilon\;\right|\;
\bnh,S_n\right]\notag\\
&=
\frac{\displaystyle\int_{W_n}^{L_n(\epsilon)}
u^{(n-p_n+a-c-2)/2}(1-u)^{(p_n+c-4)/2}\;du}
{\displaystyle\int_{W_n}^1 u^{(n-p_n+a-c-2)/2}(1-u)^{(p_n+c-4)/2}\;du}
\notag\\
&=\frac{P_M[W_n<U_n<L_n(\epsilon)\mid\bnh,S_n]}{P_M(U_n>W_n\mid\bnh
,S_n)}, \label{eq:p-g-form}
\end{align}
where $U_n\sim\text{Beta}((n-p_n+a-c)/2,\;(p_n+c-2)/2)$ and is
independent of $\bnh$ and $S_n$ under $P_M.$
Note that by the properties of the beta distribution, $U_n\to1-\alpha\,
\aspo,$ and $P_M(U_n>W_n\mid\bnh,S_n)>0$ for all~$n\,\aspo$ since
$W_n<1$ for all~$n\,\aspo$. We now introduce several technical results
regarding these quantities that will be useful in proving the main
theorem. The proofs are deferred to the Appendix.

%
\begin{lem}
\label{lem:wn-bound}
If $\liminf_{n\to\infty}\ltwo{\gn-\bon}^2\ge\delta$ for some $\delta
>0$, then $\limsup_{n\to\infty} W_n\le(1-\alpha)\lmax\sigo/(\delta+\lmax
\sigo)
<1-\alpha
\,\aspo.$
\end{lem}
%

%
\begin{lem}
\label{lem:wn-ln-to-zero}
If $\ltwo{\gn-\bon}^2\to\infty,$ then
(i)~$W_n\to0\,\aspo,$ and also
(ii)~$L_n(\epsilon)\to0\break\text{a.s.}(P_0)$ for every $\epsilon>0.$
\end{lem}
%

%
\begin{lem}
\label{lem:ln-bounds}
If $\liminf_{n\to\infty}\linf{\gn-\bon}>0$ and $\ltwo{\gn-\bon}^2\to
A$, where $0<A<\infty$, then
(i)~for every $\epsilon>0$, there exists $L^\star(\epsilon)<1$ such
that $\limsup_{n\to\infty} L_n(\epsilon)\le L^\star(\epsilon)\,\aspo$, and
(ii)~for every $\zeta\!<1,$ there exists $\epsilon_{\zeta}>0$ such that
\[
\liminf_{n\to\infty} L_n\lp\epsilon_{\zeta}\rp>\zeta\,\aspo.
\]
\end{lem}
%

To prove our main result, we will also need the following lemma, which
provides a simple result about beta random variables, the proof of
which is in the Appendix.

%
\begin{lem}
\label{lem:beta-to-mean}
Let $Z_n\sim\text{Beta}(a_n,b_n)$ for $n\ge1$, where $a_n/n\to1-\alpha$
and $b_n/n\to\alpha$, with $0\le\alpha<1$. Then $P(1-\alpha-\epsilon\le
Z_n\le1-\alpha+\epsilon)\to1$ for every $\epsilon>0$.
\end{lem}
%

We may now state and prove the main result, a necessary and sufficient
condition for posterior consistency in the hyper-$g$-prior hierarchical
model. Interestingly, this condition is identical to the one given in
Theorem~\ref{thm:pcg-eb} for the empirical Bayesian model.

%
\begin{thm}
\label{thm:pcg-conj}
In the $g$-prior model with the
hyper-$g$-prior,
posterior consistency occurs if and only if either $\alpha=0$ or there
does not exist a subsequence $k_n$ of $n$ and a constant $A>0$ such
that $||\gkn-\bokn||_2^2\to A$ and $||\gkn-\bokn||_\infty\nrightarrow0.$
\end{thm}

\noindent
The proof of this theorem is provided in the Appendix.

\subsection{Interpretations and Implications}
It should not be entirely surprising that the empirical Bayesian and
hyper-$g$-prior hierarchical models share the same necessary and
sufficient condition for posterior consistency. Indeed, the choice
$c=0$ yields the $\text{Uniform}(0,\infty)$ hyperprior on~$g,$ and in
this case the marginal posterior and likelihood of~$g$ coincide. More
generally, we should expect an adequately well-behaved hierarchical
model to exhibit broadly similar behavior to the empirical Bayesian
model, since both models essentially permit the data to determine the
value of~$g.$


\section{Zellner-Siow Hierarchical Model}

Another popular choice for the hyperprior $\pi_n(g)$ is $g\sim\text
{InverseGamma}(1/2,n/2),$ called the Zellner-Siow hyperprior \citep
{ZS1980}. The motivation behind this choice is clearest when $\xx=n\ipn
,$ in which case it leads to
marginal Cauchy priors for each component of $\bbn.$ In this section,
we will provide a sufficient condition for posterior consistency with
the Zellner-Siow hyperprior. It still remains an open problem to
determine if the condition is also necessary.

For general $\xx,$ the Zellner-Siow hyperprior yields the posterior
\begin{align}
\pi_n(g\mid\bnh,S_n)
&\propto
(g+1)^{-(p_n)/2}\lp S_n+b+\frac{T_n}{g+1}\rp^{-(n+a-2)/2}
g^{-3/2}\exp\lp-\frac{n}{2g}\rp \notag\\
&\propto
(g+1)^{(n-p_n+a-2)/2}\big[(g+1)(S_n+b)+T_n\big]^{-(n+a-2)/2}
g^{-3/2} \notag\\
&\qquad
\times\exp\lp-\frac{n}{2g}\rp. \label{eq:post-g-zs}
\end{align}
We begin with a lemma showing that Lemma~\ref{lem:pcg-hier} applies in
this model. The proof is deferred to the Appendix.

%
\begin{lem}
\label{lem:pcg-hier-zs}
With the Zellner-Siow hyperprior,
$n^{-3}\,T_n^2\,E_M[g^2(g+1)^{-4}\mid\bnh,S_n]\to0\break\text{a.s.}(P_0).$
\end{lem}
%

Now consider the form of the posterior probabilities in Lemma~\ref
{lem:pcg-hier} under this hyperprior.
By once again making the transformation in (\ref{eq:transform}), we may write
\begin{align}
&P_M\left[\left.\lp\frac{1}{g+1}\rp\linf{\gn-\bon}>\epsilon\;\right|\;
\bnh,S_n\right]
\notag\\
&\quad=
\frac{\displaystyle\int_{W_n}^{L_n(\epsilon)} \frac
{u^{(n-p_n+a-2)/2}(1-u)^{(p_n-4)/2}}{\left\{\left[u-W_n\right]/\left
[W_n(1-u)\right]\right\}^{3/2}}
\exp\left[-\frac{nW_n(1-u)}{2(u-W_n)}\right]\;du\phantom{\dfrac{\dfrac
{1}{1}}{\dfrac{1}{1}}}}
{\displaystyle\int_{W_n}^1 \frac
{u^{(n-p_n+a-2)/2}(1-u)^{(p_n-4)/2}}{\left\{\left[u-W_n\right]/\left
[W_n(1-u)\right]\right\}^{3/2}}\exp\left[-\frac
{nW_n(1-u)}{2(u-W_n)}\right]\;du\phantom{\dfrac{\dfrac{1}{1}}{\dfrac
{1}{1}}}} \notag\\
&\quad=
\frac{\displaystyle\int_{W_n}^{L_n(\epsilon)} f_n(u)\;
\left[\frac{u-W_n}{1-u}\right]^{-3/2}\exp\left[-\frac
{nW_n(1-u)}{2(u-W_n)}\right]\;du}
{\displaystyle\int_{W_n}^1 f_n(u)\;
\left[\frac{u-W_n}{1-u}\right]^{-3/2}\exp\left[-\frac
{nW_n(1-u)}{2(u-W_n)}\right]\;du},
\label{eq:g-zs-ratio}
\end{align}
where $f_n$ is the density of a $\text{Beta}[(n-p_n+a)/2,(p_n-2)/2]$
random variable with respect to Lebesgue measure.
The following lemma (proven in the Appendix) addresses the lower tail
probabilities of such a sequence.



%
\begin{lem}
\label{lem:beta-tail}
Let $Z_n\sim\text{\upshape Beta}(a_n,b_n)$ for $n\ge1,$ where $a_n/n\to
1-\alpha$ and $b_n/n\to\alpha,$ with $0\le\alpha<1$, and let $\xi\ge0$.
Then (i)~$P(Z_n\le\xi)
\le4^n\xi^{n(1-\alpha)}$ for all sufficiently large~$n$ if $\alpha>0,$
and (ii)~$P(Z_n\le\xi)\le\xi^{n/2}$ for all sufficiently large~$n$ if
$\alpha=0.$
\end{lem}
%

Note that the bound provided by Lemma~\ref{lem:beta-tail} in the case
where $0<\alpha<1$ is only useful if $\xi^{1-\alpha}<1/4.$
Now let $Q_n(\epsilon)$ and $R_n$ denote the numerator and denominator,
respectively, of (\ref{eq:g-zs-ratio}).
The following lemmas establish some results regarding these quantities
that will effectively provide the proof of the main theorem. Their
proofs are provided in the Appendix.


%
\begin{lem}
\label{lem:rn-bound}
If $\liminf_{n\to\infty}||\gn-\bon||_2^2>0,$ then there exists a finite
constant $K$ such that $R_n\ge\exp(-nK)$ for all sufficiently large~$n\,
\aspo.$
\end{lem}
%

%
\begin{lem}
\label{lem:qn-bound}
If $||\gn-\bon||_2^2\to\infty,$ then there exists a sequence of
constants $\kappa_n(\epsilon)\to\infty$ such that $Q_n(\epsilon)\le\exp
\left[-n\kappa_n(\epsilon)\right]$ for all sufficiently large~$n\,\aspo.$
\end{lem}
%

%
\begin{lem}
\label{lem:zs-alpha-zero}
If $||\gn-\bon||_2^2\to A>0,$ $\liminf_{n\to\infty}||\gn-\bon||_\infty
>0,$ and $\alpha=0,$ then $Q_n(\epsilon)/R_n\to0\,\aspo$ for every
$\epsilon>0.$
\end{lem}
%

We may now state the main theorem, which establishes the same
sufficient condition for posterior consistency under the Zellner-Siow
hyperprior as for the conjugate hyperprior and empirical Bayes models
of the previous sections. However, unlike Theorems~\ref{thm:pcg-eb}
and~\ref{thm:pcg-conj}, it does not establish the necessity of the
condition, which remains an open question.

%
\begin{thm}
\label{thm:pcg-zs}
In the $g$-prior model with the Zellner-Siow hyperprior, posterior
consistency occurs if either $\alpha=0$ or there does not exist a
subsequence $k_n$ of $n$ and a constant $A>0$ such that $||\gkn-\bokn
||_2^2\to A$ and $\linft{\gkn-\bokn}\nrightarrow0.$
\end{thm}

The proof of this theorem is provided in the Appendix.

\subsection{Interpretations and Implications}
Since the same condition is sufficient for posterior consistency under
both the hyper-$g$-prior and Zellner-Siow hierarchical models, one
might wonder if this condition is sufficient for posterior consistency
under every hierarchical model. However, the falsehood of such a claim
is made clear by the observation that the non-hierarchical model, for
which the sufficient condition differs, is simply a special case of the
hierarchical model in which the hyperprior~$\pi_n$ is specified to be
degenerate at~$g_n.$ In actuality, the posterior consistency or
inconsistency of hierarchical models with other hyperpriors on~$g$
remains a topic for future consideration.

\section{Summary}

We have derived conditions for posterior consistency under $g$-priors
by defining posterior consistency under the $\ell_\infty$ vector norm,
which allows useful results to be obtained even when the number of
parameters~$p\equiv p_n$ grows in proportion to the sample size~$n.$
Using this definition, we have obtained conditions for posterior
consistency under a variety of $g$-prior models. First, we have
obtained a necessary and sufficient condition for posterior consistency
in the non-hierarchical model in which $g\equiv g_n$ is specified as a
series of constants. Additionally, we have derived a necessary and
sufficient condition for posterior consistency under both the empirical
Bayesian $g$-prior model \citep{GF2000} and the
hyper-$g$-prior model \citep{LPMCB2008}. Interestingly, we have
found that the condition is the same for both models, and we have
illustrated that the necessity of the condition proves posterior \emph
{inconsistency} in a somewhat surprising scenario. Finally, we have
shown that this same condition is sufficient for posterior consistency
in the Zellner-Siow $g$-prior model \citep{ZS1980}, but the condition's
necessity or lack thereof remains an open question for future consideration.

\section*{Appendix: Proofs}

%
\begin{prf}
[Proof of Lemma~\ref{lem:sn}]
Under $P_0,$ the expectation and fourth central moment of $S_n$ are
$E_0(S_n)=(n-p_n)\sigo$ and $(\mu_4)_0(S_n)=12(n-p_n)(n-p_n+4)\sigma
_0^8.$ Let $\epsilon>0.$ Then
\begin{align*}
\sum_{n=1}^\infty P_0\lp\left|\frac{S_n}{n-p_n}-\sig_0\right|>\epsilon
\rp
\le
\frac{12\sigma_0^8}{\epsilon^4}\sum_{n=1}^\infty\frac
{n-p_n+4}{(n-p_n)^3}<\infty,
\end{align*}
so $(n-p_n)^{-1}S_n\to\sigo\,\aspo$ by the Borel-Cantelli lemma.
\end{prf}
%

%
\begin{prf}[Proof of Lemma~\ref{lem:tn-to-mean}]
Note that under $P_0,$ $T_n/\sigo$ has a noncentral chi-square
distribution with $p_n$~degrees of freedom and noncentrality parameter
$\frac{1}{2}n\lti||\gn-\bon||_2^2.$
Then the fourth central moment of $T_n$ under~$P_0$ is
\begin{align}
&\qquad
\lp\mu_4\rp_0\lp T_n\rp\notag\\
&\defined
E_0\left[\lp T_n-\thon\rp^4\right]
=E_0\left\{\left[T_n-\lp p_n\sigo+n\lti\ltwo{\gn-\bon}^2\rp\right
]^4\right\} \notag\\
&=
12\sigma_0^4\lp p_n\sigo+2n\lti\ltwo{\gn-\bon}^2\rp^2
+48\sigma_0^6\lp p_n\sigo+4n\lti\ltwo{\gn-\bon}^2\rp\notag\\
&\le
12\sigma_0^4\lp2p_n\sigo+2n\lti\ltwo{\gn-\bon}^2\rp^2
+48\sigma_0^6\lp4p_n\sigo+4n\lti\ltwo{\gn-\bon}^2\rp\notag\\
&=
48\sigma_0^4\thon^2+192\sigma_0^6\thon. \label{eq:nccs-fourth}
\end{align}
Define $\delta\defined\liminf_{n\to\infty}||\gn-\bon||_2^2$.
Observe that if $\alpha>0,$ then $\thon\ge p_n\sigo>\break \alpha n\sigo/2$
for all sufficiently large~$n,$ and so $\thon^{-1}=O(n^{-1}).$ If
$\delta>0,$ then $\thon\ge\break n\lti||\gn-\bon||_2^2>n\lmax^{-1}\delta/2$
for all sufficiently large~$n,$ and so $\thon^{-1}=O(n^{-1}).$ Either
way, $\thon^{-1}=O(n^{-1}),$
so the fourth central moment of $T_n/\thon$ under $P_0$ is
\begin{align*}
\lp\mu_4\rp_0\lp\frac{T_n}{\thon}\rp
&\le\frac{48\sigma_0^4}{\thon^2}+\frac{192\sigma_0^6}{\thon^3}
=O(n^{-2}).
\end{align*}
Then for any $\epsilon>0,$
\begin{align*}
\sum_{n=1}^\infty P_0\lp\left|\frac{T_n}{\thon}-1\right|>\epsilon\rp
&\le
\sum_{n=1}^\infty\frac{1}{\epsilon^4}\lp\mu_4\rp_0\lp\frac{T_n}{\thon
}\rp<\infty,
\end{align*}
which implies that $T_n/\thon\to1\,\aspo$ by the Borel-Cantelli lemma.
\end{prf}
%

%
\begin{prf}[Proof of Lemma~\ref{lem:tnt-to-mean}]
It follows from~(\ref{eq:nccs-fourth}) that the fourth central moment
of $\tnt$ under $P_0$ is
\begin{align*}
&\lp\mu_4\rp_0\lp\tnt\rp
\defined
E_0\left[\lp\tnt-\thton\rp^4\right] \\
&\qquad=
E_0\left\{\left[S_n-(n-p_n)\sigo+\frac{T_n}{g_n+1}-\frac{p_n\sigo+n\lti
\ltwo{\gn-\bon}^2}{g_n+1}\right]^4\right\} \\
&\qquad\le
8 E_0\left\{\left[S_n-E_0\lp S_n\rp\right]^4\right\}
+8 E_0\left\{\left[\frac{T_n}{g_n+1}-E_0\lp\frac{T_n}{g_n+1}\rp\right
]^4\right\} \\
&\qquad=
96(n-p_n)(n-p_n+4)\sigma_0^8
+\frac{12\sigma_0^4}{(g_n+1)^4}\lp p_n\sigo+2n\lti\ltwo{\gn-\bon}^2\rp
^2 \\
&\qquad\qquad+\frac{48\sigma_0^6}{(g_n+1)^4}\lp p_n\sigo+4n\lti\ltwo{\gn
-\bon}^2\rp\\
&\qquad\le
96(n-p_n+4)^2\sigma_0^8+48\sigma_0^4\thton^2+192\sigma_0^6\thton.
\end{align*}
Since $\thton\ge(n-p_n)\sigo,$ the fourth central moment of $\tnt/\thton
$ under $P_0$ is
\begin{align*}
&\qquad
\lp\mu_4\rp_0\lp\frac{\tnt}{\thton}\rp\\
&\le
\frac{96(n-p_n+4)^2\sigma_0^8}{\thton^4}+\frac{48\sigma_0^4}{\thton
^2}+\frac{192\sigma_0^6}{\thton^3}
\le
\frac{96(n-p_n+4)^2}{(n-p_n)^4}+\frac{48}{(n-p_n)^2}+\frac{192}{(n-p_n)^3},
\end{align*}
which is $O(n^{-2}).$ Then for any $\epsilon>0,$
\begin{align*}
\sum_{n=1}^\infty P_0\lp\left|\frac{\tnt}{\thton}-1\right|>\epsilon\rp
&\le
\sum_{n=1}^\infty\frac{1}{\epsilon^4}\lp\mu_4\rp_0\lp\frac{\tnt}{\thton
}\rp<\infty,
\end{align*}
which implies that $\tnt/\thton\to1\,\aspo$ by the Borel-Cantelli lemma.
\end{prf}
%

%
\begin{prf}[Proof of Lemma~\ref{lem:normal-means-noncvgc}]
Assume $\linf{\mun-\xin}\nrightarrow0.$ Then there exists a subsequence
$k_n$ of $n$ and a $\delta>0$ such that $\linf{\mukn-\xikn}>\delta$ for
all $n.$ There also exists an $i_n,$ $1\le i_n\le p_n,$ such that $\left
|\mu_{k_n,i_n}-\theta_{k_n,i_n}\right|=\linf{\mukn-\xikn}>\delta$ for
all $n.$\vadjust{\goodbreak}
Then either $\mu_{k_n,i_n}<\theta_{k_n,i_n}-\delta$ (Case~1) or $\mu
_{k_n,i_n}>\theta_{k_n,i_n}+\delta$ (Case~2).
Now let $0<\epsilon<\delta,$ and note that
$
P\lp\linf{\zkn-\xikn}\le\epsilon\rp
\le
P\lp\xi_{k_n,i_n}-\delta\le Z_{k_n,i_n}\le\xi_{k_n,i_n}+\delta\rp.
$
Recall that $\Sigmakn$ is assumed positive definite.
Then in Case~1,
\[
P\lp\xi_{k_n,i_n}-\delta\le Z_{k_n,i_n}\le\xi_{k_n,i_n}+\delta\rp
\le P\lp\mu_{k_n,i_n}\le Z_{k_n,i_n}\rp=1/2,
\]
while in Case~2,
\[
P\lp\xi_{k_n,i_n}-\delta\le Z_{k_n,i_n}\le\xi_{k_n,i_n}+\delta\rp
\le P\lp Z_{k_n,i_n}\le\mu_{k_n,i_n}\rp=1/2.
\]
%
Either way, $P(\linf{\zkn-\xikn}>\epsilon)\ge1/2$ for all $n.$
\end{prf}
%

%
\begin{prf}[Proof of Lemma~\ref{lem:normal-center}]
Note that for any $t>0,$ $\Phi(z+t)-\Phi(z-t)$ is maximized at $z=0.$ Hence,
\begin{align*}
P(|Z|\le\xi)=P(-\xi\le Z\le\xi)
&= P\lp\frac{-\xi-\mu}{\tau}\le\frac{Z-\mu}{\tau}\le\frac{\xi-\mu}{\tau
}\rp\\
&=\Phi\lp\frac{\xi}{\tau}-\frac{\mu}{\tau}\rp-\Phi\lp-\frac{\xi}{\tau
}-\frac{\mu}{\tau}\rp
\\
&\le\Phi\lp\frac{\xi}{\tau}\rp-\Phi\lp-\frac{\xi}{\tau}\rp,
\end{align*}
from which it immediately follows that $P(|Z|\le\xi)\le1-2\Phi(-\xi/\tau).$
\end{prf}
%

%
\begin{prf}[Proof of Lemma~\ref{lem:normal-correlations}]
For each $i=1,\ldots,p_n,$ partition $\bm{\Sigma}_n$ as
\[
\bm{\Sigma}_n=
\begin{bmatrix}
\bm{\Sigma}_{n,i,11} & \bm{\Sigma}_{n,i,1i} & \bm{\Sigma}_{n,i,12} \\
\bm{\Sigma}_{n,i,1i}^T & \Sigma_{n,ii} & \bm{\Sigma}_{n,i,2i}\\
\bm{\Sigma}_{n,i,12}^T & \bm{\Sigma}_{n,i,2i}^T & \bm{\Sigma}_{n,i,22}
\end{bmatrix}
,
\]
where the submatrices $\bm{\Sigma}_{n,i,11}$ and $\bm{\Sigma}_{n,i,22}$
along the diagonal have dimension $(i-1)\times(i-1)$ and $(p_n-i)\times
(p_n-i),$ respectively. Then define
$
\widetilde{\Sigma}_{n,i}\defined\text{Var}(Z_i\mid Z_{i+1},\ldots, Z_{p_n}),
$
so that
$
\widetilde{\Sigma}_{n,i}=\Sigma_{n,ii}-\bm{\Sigma}_{n,i,2i}\bm{\Sigma
}_{n,i,22}^{-1}\bm{\Sigma}_{n,i,2i}^T.
$
Note that $\widetilde{\Sigma}_{n,i}^{-1}$ is the first diagonal entry of
\[
\begin{bmatrix}
\Sigma_{n,ii} & \bm{\Sigma}_{n,i,2i}\\
\bm{\Sigma}_{n,i,2i}^T & \bm{\Sigma}_{n,i,22}
\end{bmatrix}
^{-1},
\]
which has eigenvalues bounded above by $\omega_{\min}^{-1}$ since the
eigenvalues of a principal submatrix are bounded below by the smallest
eigenvalue of the full matrix. Hence $\widetilde{\Sigma}_{n,i}^{-1}\le
\omega_{\min}^{-1},$ and the result immediately follows.
\end{prf}
%

%
\begin{prf}[Proof of Lemma~\ref{lem:sig-bounds}]
Recall that $\tnt/\thton\to1\,\aspo$ by Lemma~\ref{lem:tnt-to-mean}.
Then for all sufficiently large~$n$,
%
\begin{align*}
&\qquad P_M\lp\left. \frac{\thton}{2n}\le\sig\le\frac{2\thton}{n}\given
\bnh,S_n\rp\\
&\ge
P_M\lp\left. \frac{3\tnt}{4(n+a-4)}\le\sig\le\frac{5\tnt
}{4(n+a-4)}\given\bnh,S_n\rp
\aspo\\
&=
P_M\lp\left. \left|\sig-E_M\lp\sig\mid\bnh,S_n\rp\right|\le\frac{\tnt
}{4(n+a-4)}\given\bnh,S_n\rp\\
&\ge
1-\lp\frac{4(n+a-4)}{\tnt}\rp^2\lp\frac{2\tnt^2}{(n+a-4)^2(n+a-6)}\rp\\
&=
1-\frac{32}{n+a-6}\to1,
\end{align*}
where the last inequality is a consequence of Chebyshev's inequality,
for which we note that $\text{Var}_M(\sig\mid\bnh
,S_n)=2(n+a-4)^{-2}(n+a-6)^{-1}\tnt^2.$
\end{prf}
%

%
\begin{prf}[Proof of Theorem~\ref{thm:fixed-g}]
By Lemma~\ref{lem:pcg-basic}, we may replace $\bon$ with $\bnh$ in the
definition of posterior consistency. We will now consider four cases.

Case 1: Suppose $(g_n+1)^{-1}||\gn-\bon||_\infty\nrightarrow0.$ Then
since $||\gn-\bnh||_\infty\ge||\gn-\bon||_\infty-||\bnh-\bon||_\infty$
and $||\bnh-\bon||_\infty\to0\,\aspo$ by Lemma~\ref{lem:mle-linf}, it
follows that
$(g_n+1)^{-1}||\gn-\bnh||_\infty\nrightarrow0\,\aspo.$
Now observe that under~$P_M,$
\[
\bbn-\bnh\mid\sig,\bnh,S_n\sim N_{p_n}\lp\frac{1}{g_n+1}\lp\gn-\bnh\rp
,\frac{g_n}{g_n+1}\sig\xxi\rp.
\]
Then by Lemma~\ref{lem:normal-means-noncvgc}, there exists an $\epsilon
>0$ and a subsequence $k_n$ of $n$ such that, $\aspo$, $P_M(||\bm{\beta
}_{k_n}-\bknh||_\infty>\epsilon\mid\sig,\bknh,S_{k_n})>1/2$ for every
$n$ and every $\sig>0.$ Then
\begin{align*}
&\qquad P_M\lp||\bm{\beta}_{k_n}-\bknh||_\infty>\epsilon\mid\bknh
,S_{k_n}\rp\\
&=
E_M\left[\left.P_M\lp||\bm{\beta}_{k_n}-\bknh||_\infty>\epsilon\mid\sig
,\bknh,S_{k_n}\rp\;\right|\;\bknh,S_{k_n}\right] \\
&\ge1/2\quad\text{ for every $n$}\; \aspo.
\end{align*}
Therefore $P_M(||\bbn-\bnh||_\infty>\epsilon\mid\bnh,S_n)\nrightarrow
0,$ so posterior consistency does not occur.

For the remaining cases, suppose $(g_n+1)^{-1}||\gn-\bon||_\infty\to0.$
Then since $||\gn-\bnh||_\infty\le||\gn-\bon||_\infty+||\bnh-\bon
||_\infty$ and $||\bnh-\bon||_\infty\to0\,\aspo$ by Lemma~\ref
{lem:mle-linf}, it follows that
$(g_n+1)^{-1}||\gn-\bnh||_\infty\to0\,\aspo.$ Then
\begin{align*}
&P_M\lp||\bbn-\bnhb||_\infty>2\epsilon\mid\bnh,S_n\rp
-P_M\lp||\bnhb-\bnh||_\infty>\epsilon\mid\bnh,S_n\rp\\
&\qquad\le
P_M\lp||\bbn-\bnh||_\infty>\epsilon\mid\bnh,S_n\rp\\
&\qquad\qquad\le
P_M\lp||\bbn-\bnhb||_\infty>\epsilon/2\mid\bnh,S_n\rp+P_M\lp||\bnhb-\bnh
||_\infty>\epsilon/2\mid\bnh,S_n\rp
\end{align*}
by the triangle inequality. Note that
\[
P_M\lp||\bnhb-\bnh||_\infty>\epsilon\mid\bnh,S_n\rp=I(||\bnhb-\bnh
||_\infty>\epsilon),
\]
where $I(\cdot)$ denotes the indicator function. But $\bnhb-\bnh
=(g_n+1)^{-1}(\gn-\bnh),$ so this indicator is zero for all
sufficiently large $n\,\aspo.$ Therefore, posterior consistency occurs
in Cases~2--3 below if and only if $P_M(||\bbn-\bnhb||_\infty>\epsilon
\mid\bnh,S_n)\to0\,\aspo$ for every $\epsilon>0.$ We now consider the
individual cases.\vadjust{\goodbreak}

Case 2: Suppose that $(g_n+1)^{-1}||\gn-\bon||_\infty\to0,$ and also
suppose that $g_n(g_n+1)^{-2}(\log p_n)n^{-1}||\gn-\bon||_2^2\to0.$
Observe that
\begin{align*}
&P_M\lp||\bbn-\bnhb||_\infty>\epsilon\mid\bnh,S_n\rp\notag\\
&\qquad=
E_M\left[\left.P_M\lp\linf{\bbn-\bnhb}>\epsilon\mid\sig,\bnh,S_n\rp\;
\right|\;\bnh,S_n\right] \notag\\
&\qquad\le
E_M\left[\left.P_M\lp\linf{\bbn-\bnhb}>\epsilon\mid\sig,\bnh,S_n\rp I\lp
\sig\le\frac{2\thton}{n}\rp\;\right|\;\bnh,S_n\right]\notag\\
&\qquad\qquad+ P_M\lp\left.\sig>\frac{2\thton}{n}\given\bnh,S_n\rp
\notag.
\end{align*}
We immediately have that $P_M(\sig>2\thton/n\mid\bnh,S_n)\to0\,\aspo$
by Lemma~\ref{lem:sig-bounds}, so it suffices to work with the first
term to establish posterior consistency.
Let $v_{n,ij}$ denote the $ij$th element of $\nxxismall,$ and note
specifically that the diagonal elements may be bounded by $\lmin\le
v_{n,ii}\le\lmax$ for all~$n$ and~$i.$ Also recall that $\bbn-\bnhb\mid
\sig,\bnh,S_n\sim N_{p_n}(\bm{0}_{p_n},g_n(g_n+1)^{-1}\sig\xxit)$ under
$P_M.$ Now let $\epsilon>0,$ and bound the aforementioned first term by
\begin{align*}
& E_M\left[\left.P_M\lp\linf{\bbn-\bnhb}>\epsilon\mid\sig,\bnh,S_n\rp
I\lp\sig\le\frac{2\thton}{n}\rp\;\right|\;\bnh,S_n\right] \\
&\qquad\le
E_M\left[\left.\sum_{i=1}^{p_n}P_M\lp\left|\beta_{n,i}-\hat{\beta
}_{n,i}^{\text{B}}\right|>\epsilon\mid\sig,\bnh,S_n\rp I\lp\sig\le\frac
{2\thton}{n}\rp\;\right|\;\bnh,S_n\right] \\
&\qquad\le
E_M\left[\left.\sum_{i=1}^{p_n}2\Phi\lp-\sqrt{\frac{\epsilon^2
(g_n+1)n}{g_n v_{n,ii}\sig}}\rp I\lp\sig\le\frac{2\thton}{n}\rp\;\right
|\;\bnh,S_n\right] \\
&\qquad\le
2p_n E_M\left[\left.\Phi\lp-\sqrt{\frac{\epsilon^2 (g_n+1)n^2}{2g_n\lmax
\thton}}\rp\;\right|\;\bnh,S_n\right]
=
2p_n \Phi\lp-\sqrt{\frac{\epsilon^2 (g_n+1)n^2}{2g_n\lmax\thton}}\rp
\end{align*}
where $\Phi(\cdot)$ denotes the standard normal cdf.
Then by the Mills ratio,
\begin{align*}
2p_n \Phi\lp-\sqrt{\frac{\epsilon^2 (g_n+1)n^2}{2g_n\lmax\thton}}\rp
&\le
2p_n\sqrt{\frac{g_n\lmax\thton}{\pi\epsilon^2 (g_n+1)n^2}}\;\exp\lp
-\frac{\epsilon^2 (g_n+1)n^2}{4g_n\lmax\thton}\rp.
\end{align*}
This expression clearly tends to zero if $\thton/n$ is bounded above,
so we may instead assume that $\thton/n\to\infty,$ which by inspection
occurs if and only if $(g_n+1)^{-1}||\gn-\bon||_2^2\to\infty.$ Then
$\thton\le2n\lti(g_n+1)^{-1}||\gn-\bon||_2^2$ for all sufficiently
large $n,$ and hence
\begin{align*}
&\qquad2p_n \Phi\lp-\sqrt{\frac{\epsilon^2 (g_n+1)n^2}{2g_n\lmax\thton
}}\rp\\
&\le
2p_n\sqrt{\frac{2\lmax\; g_n\ltwo{\gn-\bon}^2}{\pi\lton\epsilon^2
(g_n+1)^2\,n}}\;\exp\lp-\frac{\lton\epsilon^2 (g_n+1)^2\,n}{8\lmax\;
g_n\ltwo{\gn-\bon}^2}\rp\\
&\le
\sqrt{\frac{8\lmax\; g_n\ltwo{\gn-\bon}^2\;\log p_n}{\pi\lton\epsilon^2
(g_n+1)^2\,n}}\\
&\qquad\times\exp\left[\lp1-\frac{\lton\epsilon^2 (g_n+1)^2\,n}{8\lmax
\; g_n\ltwo{\gn-\bon}^2\;\log p_n}\rp\log p_n\right]\\
&\to0\;\text{ for every $\epsilon>0$}
\end{align*}
by the assumption that $g_n(g_n+1)^{-2}(\log p_n)n^{-1}||\gn-\bon
||_2^2\to0.$ Therefore, posterior consistency occurs.

Case 3: Suppose $(g_n+1)^{-1}||\gn-\bon||_\infty\to0,$ but now suppose
that $g_n(g_n+1)^{-2}(\log p_n)n^{-1}||\gn-\bon||_2^2\nrightarrow0.$
Then there exist a subsequence $k_n$ of $n$ and a constant $\delta>0$
such that $g_{k_n}(g_{k_n}+1)^{-2}(\log p_{k_n})k_n^{-1}||\gkn-\bokn
||_2^2>\delta$ for all~$n.$ Note that posterior inconsistency of the
subsequence $P_M(\bbkn\mid\bknh,S_{k_n})$ implies posterior
inconsistency of the overall sequence $P_M(\bbn\mid\bnh,S_{n}),$ so we
may assume without loss of generality that $k_n=n$ for notational convenience.
Also, define $\Sigman$ to be the $p_{n}\times p_{n}$ matrix with
elements $\Sigma_{n,ij}\defined v_{n,ij}/\sqrt{v_{n,ii}v_{n,jj}},$
where $v_{n,ij}$ denotes the $ij$th element of $\nxxismall$ as before. Then
\begin{align*}
& P_M\lp||\bbn-\bnhb||_\infty>\epsilon\mid\bnh,S_n\rp\\
&\ge
E_M\left[\left.P_M\lp\linf{\bbn-\bnhb}>\epsilon\mid\sig,\bnh,S_n\rp I\lp
\sig\ge\frac{\thton}{2n}\rp\;\right|\;\bnh,S_n\right]\\
&\ge
E_M\Bigg[\left.P_M\lp\left.\max_{1\le i\le p_{n}}\left|\beta_{n,i}-\hat
{\beta}_{n,i}^\text{B}\right|>\sqrt{\epsilon^2\frac{v_{{n},i}}{\lmin
}}\given\sig,\bnh,S_{n}\rp
I\lp\sig\ge\frac{\thton}{2n}\rp\;\right|\;\bnh,S_{n}\Bigg].
\end{align*}
Then we may write
\begin{align*}
& P_M\lp||\bbn-\bnhb||_\infty>\epsilon\mid\bnh,S_n\rp\\
&\qquad\ge
E_M\left[\left.P_M\lp\left.\max_{1\le i\le p_{n}}\left|Z_i\right|>\sqrt
{\frac{(g_{n}+1)\epsilon^2 n}{g_{n}\lmin\sig}}\given\sig
\rp
I\lp\sig\ge\frac{\thton}{2n}\rp\;\right|\;\bnh,S_{n}\right],
\end{align*}
where $\bm{Z}_{n}
\sim N_{p_{n}}(\bm{0}_{p_{n}},\bm{\Sigma}_{n})$ and is independent of
$\sig$ under $P_M.$ Now note that the innermost conditional probability
is a nondecreasing function of $\sig$, which implies that
%
\begin{align*}
&P_M\lp||\bbn-\bnhb||_\infty>\epsilon\mid\bnh,S_{n}\rp\\
&\qquad\ge
E_M\left[\left.P_M\lp\left.\max_{1\le i\le p_{n}}\left|Z_i\right|>\sqrt
{\frac{2(g_{n}+1)\epsilon^2 n^2}{g_{n}\lmin\thton}}\given\sig
\rp I\lp\sig\ge\frac{\thton}{2n}\rp\;\right|\;\bnh,S_{n}\right] \\
&\qquad=
P_M\lp\max_{1\le i\le p_{n}}\left|Z_i\right|>\sqrt{\frac
{2(g_{n}+1)\epsilon^2 n^2}{g_{n}\lmin\thton}}\rp P_M\lp\left.\sig\ge
\frac{\thton}{2n}\given\bnh,S_n\rp,
\end{align*}
since the entries of $\bm{\Sigma}_n$ depend only on $\xx.$
Then Lemma~\ref{lem:sig-bounds} immediately implies that $P_M(\sig\ge
\thton/2n\mid\bnh,S_n)\to1\,\aspo,$ so it suffices to show that the
first term is bounded away from zero for all sufficiently large $n$.
Now define $\etaon\defined[2(g_n+1)\epsilon^2 n^2/g_n\lmin\thton
]^{1/2}$ and $\widetilde{\Sigma}_{n,i}\defined\text{Var}(Z_i\mid
Z_{i+1},\ldots,Z_{p_n}).$
Then
%
\begin{align*}
P_M\lp\max_{1\le i\le p_n}|Z_i|\le\etaon\rp
&=
E_M\!\left[P_M\lp\left.|Z_1|\le\etaon\given Z_2,Z_3,\ldots,Z_{p_n}\rp
\prod_{i=2}^{p_n}I_{\{|Z_i|\le\etaon\}}\right] \\
&\le
\left[1-2\Phi\lp-\etaon/\sqrt{\widetilde{\Sigma}_{n,1}}\rp\right]
E_M\!\left[
\prod_{i=2}^{p_n}I_{\{|Z_i|\le\etaon\}}\right]
\end{align*}
%
by Lemma~\ref{lem:normal-center} and the fact that $\widetilde{\Sigma
}_{n,1}$ does not depend on $Z_2,Z_3,\ldots,Z_{p_n}.$
By repeated conditioning on $Z_{i+1},Z_{i+2},\ldots,Z_{p_n}$ for
$i=2,3,\ldots,p_n-1$ and application of Lemma~\ref{lem:normal-center}
as above, we find that
%
\begin{align*}
P_M\lp\max_{1\le i\le p_n}|Z_i|\le\etaon\rp
\le
\prod_{i=1}^{p_n}\left[1-2\Phi\lp-\etaon/\sqrt{\widetilde{\Sigma
}_{n,i}}\rp\right].
\end{align*}
%
Note that
\begin{align*}
\thton\ge
\frac{n\lti}{g_n+1}\ltwo{\gn-\bon}^2
\ge
\frac{\delta(g_n+1)n^2}{\lmax\,g_n\log p_n},
\end{align*}
which implies that
\begin{align*}
\etaon\le\sqrt{\frac{2\lmax\epsilon^2\log p_n}{\delta\lmin}}.
\end{align*}
The eigenvalues of $\bm{\Sigma}_n$ are bounded below by $\lmin/\lmax,$
so $\inf_{n,i}\widetilde{\Sigma}_{n,i}\ge\lmin/\lmax$ by Lemma~\ref
{lem:normal-correlations}. Then it follows that
%
\begin{align*}
P_M\lp\max_{1\le i\le p_n}|Z_i|\le\etaon\rp
&\le
\left[1-2\Phi\lp-\sqrt{\frac{2\lmax^2\epsilon^2\log p_n}{\delta\lmin
^2}}\rp\right]^{p_n}\\
&\le
\exp\!\left[-2p_n\Phi\lp-\sqrt{\frac{2\lmax^2\epsilon^2\log p_n}{\delta
\lmin^2}}\rp\right].
\end{align*}
%
Notice that if any subsequence of $p_n$ is bounded above, then the quantity
\[
\Phi[-(2\lmax^2\epsilon^2\log p_n/\delta\lmin^2)^{1/2}]
\]

\noindent
is bounded away from zero along that subsequence, and thus posterior
inconsistency follows immediately. So we may instead assume that $p_n\to
\infty.$ Then
\[
2\lmax^2\epsilon^2\log p_n/\delta\lmin^2\to\infty,
\]

\noindent
in which case the inequality
\[
1-\Phi(t)\ge(t^{-1}-t^{-3})(2\pi)^{-1/2}\exp(-t^2/2)
\ge(2t)^{-1}(2\pi)^{-1/2}\exp(-t^2/2)
\]
for large~$t$ may be applied for all sufficiently large $n,$ yielding
%
\begin{align*}
P_M\lp\max_{1\le i\le p_{n}}\left|Z_i\right|>\etaon\rp
&\ge
1-
\exp\left[-p_n\;
\sqrt{\frac{\delta\lmin^2}{4\pi\lmax^2\epsilon^2\log p_n}}\;
\exp\lp-\frac{\lmax^2\epsilon^2\log p_n}{\delta\lmin^2}\rp\right] \\
&=
1-
\exp\left\{-
\sqrt{\frac{\delta\lmin^2}{4\pi\lmax^2\epsilon^2\log p_n}}\;
\exp\left[\lp1-\frac{\lmax^2\epsilon^2}{\delta\lmin^2}\rp\log p_n\right
]\right\}\\
&\to1\;\text{ for $\epsilon<\sqrt{\frac{\delta\lmin^2}{\lmax^2}}$}.
\end{align*}
%
Therefore posterior consistency does not occur.
\end{prf}

%
\begin{prf}[Proof of Lemma~\ref{lem:eb-g-bound}]
Define $\delta\defined\liminf_{n\to\infty}||\gn-\bon||_2^2,$ and assume
$\delta>0.$ Then
\[
\frac{T_n}{p_n}=\frac{T_n}{\thon}\lp\sigo+\frac{n}{p_n\lton}\ltwo{\gn
-\bon}^2\rp
>
\frac{T_n}{\thon}\lp\sigo+\frac{\delta}{2\lmax}\rp> \sigo+\frac{\delta
}{4\lmax}
\]
for all sufficiently large~$n\,\aspo,$ since $T_n/\thon\to1\,\aspo$ by
Lemma~\ref{lem:tn-to-mean}. Then
\begin{align*}
\liminf_{n\to\infty}\gnheb&>\liminf_{n\to\infty}\left[\lp\frac
{n-p_n+a-2}{S_n+b}\rp\lp\sigo+\frac{\delta}{4\lmax}\rp-1\right]
=\frac{\delta}{4\lmax\sigo}>0\;\;\aspo
\end{align*}
since $(n-p_n+a-2)/(S_n+b)\to1/\sigo\,\aspo$ by Lemma~\ref{lem:sn}.
\end{prf}
%

%
\begin{prf}[Proof of Theorem~\ref{thm:pcg-eb}]
By Theorem~\ref{thm:fixed-g}, we immediately have that posterior
consistency occurs if and only if both
\begin{equation}
\frac{||\gn-\bon||_\infty}{\gnheb+1}\to0\;\text{ and }\;\frac{\gnheb
\log p_n}{(\gnheb+1)^2 \,n}||\gn-\bon||_2^2\to0\;\aspo. \label{eq:eb-basic}
\end{equation}
We now consider three cases.

Case~1: Suppose there do not exist a subsequence $k_n$ of $n$ and a
constant $A>0$ such that $||\gkn-\bokn||_2^2\to A$ and $||\gkn-\bokn
||_\infty\nrightarrow0.$ Now let $k_n$ be a subsequence of $n,$ and
consider two sub-cases.

Case~1.1: Suppose $||\gkn-\bon||_\infty\to0.$ Then clearly the first
condition in~(\ref{eq:eb-basic}) is satisfied trivially. Note that for
any further subsequence~$m_n$ of~$k_n$ for which $||\gmn-\bomn||_2^2\to
0,$ the second condition in~(\ref{eq:eb-basic}) is satisfied trivially
as well, so we may instead assume $\liminf_{n\to\infty}||\gkn-\bokn
||_2^2>0.$ Then for all sufficiently large~$n\,\aspo,$
\begin{align}
&\qquad \frac{\gknheb\log p_{k_n}}{(\gknheb+1)^2 \,k_n}||\gkn-\bokn
||_2^2 \notag\\
&\le
\frac{\log k_n}{k_n}\lp\frac{S_{k_n}+b}{k_n-p_{k_n}+a-2}\rp\lp\frac
{\thokn}{T_{k_n}}\rp
\lp\frac{p_{k_n}||\gkn-\bokn||_2^2}{\thokn}\rp\notag\\
&\le
\frac{\log k_n}{k_n}\lp\frac{S_{k_n}+b}{k_n-p_{k_n}+a-2}\rp\lp\frac
{\thokn}{T_{k_n}}\rp
\frac{p_{k_n}\lmax}{k_n}\to0\;\;\aspo \label{eq:eb-second}
\end{align}
by Lemmas~\ref{lem:sn}, \ref{lem:tn-to-mean}, and~\ref{lem:eb-g-bound}.
Thus, both conditions in~(\ref{eq:eb-basic}) hold along the
subsequence~$k_n$.\vadjust{\goodbreak}

Case~1.2: Note that Case~1.1 can be applied to any further
subsequence~$m_n$ of $k_n$ for which $||\gmn-\bomn||_\infty\to0,$ so we
may suppose for Case~1.2 that $\liminf_{n\to\infty}||\gkn-\bokn||_\infty
>0.$ Note also that in this case, there cannot exist any further
subsequence~$m_n$ of~$k_n$ for which $||\gmn-\bomn||_2^2$ converges to
a nonzero constant, since this would contradict the original
supposition of Case~1. Then since $\liminf_{n\to\infty}||\gkn-\bokn
||_2^2\ge\liminf_{n\to\infty}||\gkn-\bokn||_\infty^2>0,$ it follows
that $||\gkn-\bokn||_2^2\to\infty.$ Then for all sufficiently large~$n\,
\aspo,$
\begin{align}
\label{eq:eb-first}
&\qquad\frac{1}{\gknheb+1}||\gkn-\bokn||_\infty\notag\\
&=
\lp\frac{S_{k_n}+b}{k_n-p_{k_n}+a-2}\rp
\lp\frac{\thokn}{T_{k_n}}\rp\lp\frac{p_{k_n}\linf{\gkn-\bokn}}{\thokn
}\rp\\
&\le
\lp\frac{S_{k_n}+b}{k_n-p_{k_n}+a-2}\rp
\lp\frac{\thokn}{T_{k_n}}\rp\lp\frac{p_{k_n}\lmax}{k_n\ltwo{\gkn-\bokn
}}\rp
\to0\;\aspo \notag
\end{align}
by Lemmas~\ref{lem:sn}, \ref{lem:tn-to-mean}, and~\ref{lem:eb-g-bound},
while (\ref{eq:eb-second}) also holds by the same lemmas. Thus, both
conditions hold along the subsequence~$k_n.$ Since Cases~1.1 and~1.2
together establish that both conditions hold along any
subsequence~$k_n$, they hold for the whole sequence, and therefore
posterior consistency occurs.

Case 2: Now suppose there exist a subsequence $k_n$ of $n$ and a
constant $A>0$ such that $||\gkn-\bokn||_2^2\to A>0$ and $||\gkn-\bokn
||_\infty\nrightarrow0,$ and suppose $\alpha=0.$
Note that Case~1.1 can be applied to any further subsequence~$m_n$ of
$k_n$ for which $||\gmn-\bomn||_\infty\to0,$ so we may suppose for
Case~2 that $\liminf_{n\to\infty}||\gkn-\bokn||_\infty>0.$ Then (\ref
{eq:eb-second}) and (\ref{eq:eb-first}) still hold by Lemmas~\ref
{lem:sn}, \ref{lem:tn-to-mean}, and~\ref{lem:eb-g-bound} since
$p_{k_n}/k_n\to0$ and $\liminf_{n\to\infty}||\gkn-\bokn||_2\ge\liminf
_{n\to\infty}||\gkn-\bokn||_\infty>0.$
Hence, the two conditions hold for every subsequence, and consequently
for the overall sequence. Therefore posterior consistency
occurs.

Case 3: Now suppose there exist a subsequence $k_n$ of $n$ and a
constant $A>0$ such that $||\gkn-\bokn||_2^2\to A>0$ and $||\gkn-\bokn
||_\infty\nrightarrow0,$ but suppose $\alpha>0.$
As in Case~2, we may suppose for Case~3 that $\liminf_{n\to\infty}||\gkn
-\bokn||_\infty>0.$
Then for all sufficiently large~$n\,\aspo,$
\begin{align*}
\frac{||\gkn-\bokn||_\infty}{\gknheb+1}
&=
\lp\frac{S_{k_n}+b}{k_n-p_{k_n}+a-2}\rp
\lp\frac{\thokn}{T_{k_n}}\rp\lp\frac{p_{k_n}\linf{\gkn-\bokn}}{\thokn
}\rp\\
&\ge
\lp\frac{S_{k_n}+b}{k_n-p_{k_n}+a-2}\rp
\lp\frac{\thokn}{T_{k_n}}\rp\lp\frac{p_{k_n}\lmin}{k_n\ltwo{\gkn-\bokn
}^2}\rp\\
&\qquad\quad\times\liminf_{n\to\infty}\linf{\gkn-\bokn} \\
&\to\frac{\sigo\alpha\lmin}{A}\liminf_{n\to\infty}||\gkn-\bokn||_\infty>0
\quad\aspo
\end{align*}
by Lemmas~\ref{lem:sn}, \ref{lem:tn-to-mean}, and~\ref{lem:eb-g-bound}.
The first condition fails for the subsequence $k_n$ and hence for the
overall sequence. Therefore posterior consistency does not
occur.
\end{prf}

%
\begin{prf}[Proof of Lemma~\ref{lem:pcg-hier}]
Assume that $n^{-3}\,T_n^2\,E_M[g^2(g+1)^{-4}\mid\bnh,S_n]\to0\,\aspo.$
By Lemma~\ref{lem:pcg-basic}, to determine whether posterior
consistency occurs, it suffices to consider whether $P_M(||\bbn-\bnh
||_\infty>\epsilon\mid\bnh,S_n)\to0\,\aspo$ for every $\epsilon>0.$
By iterated expectation and the triangle inequality,
\begin{align}
&E_M\left[\left.P_M\lp\left.\linf{\bntbg-\bnh}>2\epsilon\given g,\sig
,\bnh,S_n\rp\given\bnh,S_n\right] \notag\\
&-E_M\left[\left.P_M\lp\left.\linf{\bbn-\bntbg}>\epsilon\given g,\sig
,\bnh,S_n\rp\given\bnh,S_n\right]\notag\\
&\quad\le
E_M\left[\left.P_M\lp\left.\linf{\bbn-\bnh}>\epsilon\given g,\sig,\bnh
,S_n\rp\given\bnh,S_n\right]\label{eq:tri-ineq-conj-1}\\
&\quad\le
E_M\left[\left.P_M\lp\left.\linf{\bntbg-\bnh}>\epsilon/2\given g,\sig
,\bnh,S_n\rp\given\bnh,S_n\right] \notag\\
&\quad\quad+E_M\left[\left.P_M\lp\left.\linf{\bbn-\bntbg}>\epsilon
/2\given g,\sig,\bnh,S_n\rp\given\bnh,S_n\right]. \notag
\end{align}
Consider $P_M(||\bbn-\bntbg||_\infty>\epsilon\mid g,\sig,\bnh,S_n)$ for
some arbitrary $\epsilon>0$ and $g\ge0.$ Under $P_M,$
\begin{align*}
\left.\bbn-\bntbg\given g,\sig,\bnh,S_n\sim N_{p_n}\lp\bm{0},\frac
{g}{g+1}\sig\xxi\rp.
\end{align*}
Let $v_{n,11},\ldots,v_{n,p_np_n}$ denote the diagonal elements of
$\nxxismall,$ and write
\begin{align*}
&\qquad P_M\lp\left.\linf{\bbn-\bntbg}>\epsilon\given g,\sig,\bnh,S_n\rp
\\
&\leq
\sum_{i=1}^{p_n}P_M\lp\left.\left|\beta_{n,i}-\tilde{\beta}_{n,i}^\text
{B}(g)\right|>\epsilon\given g,\sig,\bnh,S_n\rp\\
&\leq
\sum_{i=1}^{p_n}P_M\lp\left.\left[\beta_{n,i}-\tilde{\beta}_{n,i}^\text
{B}(g)\right]^4>\epsilon^4\given g,\sig,\bnh,S_n\rp\\
&\leq
\sum_{i=1}^{p_n}\frac{3g^2\sigma^4 v_{n,ii}}{(g+1)^2 n^2\epsilon^4}
\leq
\frac{3\lmax g^2\sigma^4}{(g+1)^2 n\epsilon^4}.
\end{align*}
Then
\begin{align*}
&E_M\left[\left.P_M\lp\left.\linf{\bbn-\bntbg}>\epsilon\given g,\sig
,\bnh,S_n\rp\given\bnh,S_n\right] \\
&\qquad\le
\frac{3\lmax}{n\epsilon^4}E_M\lp\left.\frac{g^2\sigma^4}{(g+1)^2}\given
\bnh,S_n\rp\\
&\qquad=
\frac{3\lmax}{n\epsilon^4}E_M\left[\left.\frac{g^2}{(g+1)^2}E_M\lp\left
.\sigma^4\given g,\bnh,S_n\rp\given\bnh,S_n\right].
\end{align*}
Observe from the form of the posterior in (\ref{eq:post-sg}) that under $P_M,$
\[
\sig\mid g,\bnh,S_n\sim\text{InverseGamma}\lp\frac{n+a-2}{2},\frac
{S_n+b+(g_n+1)^{-1}T_n}{2}\rp.\vadjust{\eject}
\]
Therefore,
\begin{align*}
&E_M\left[\left.P_M\lp\left.\linf{\bbn-\bntbg}>\epsilon\given g,\sig
,\bnh,S_n\rp\given\bnh,S_n\right] \\
&=
\frac{3\lmax}{n\epsilon^4}E_M\left[\left.\frac{g^2\left[
S_n+b+(g_n+1)^{-1}T_n\right]^2}{(g+1)^2(n+a-4)(n+a-6)}\given\bnh
,S_n\right] \\
&\le
\frac{6\lmax}{n\epsilon^4}\lp\frac{S_n+b}{n+a-6}\rp^2
+\frac{6\lmax}{n\epsilon^4}\frac{T_n^2}{(n+a-6)^2}E_M\left[\left.\frac
{g^2}{(g+1)^4}\given\bnh,S_n\right]\to0\;\;\aspo
\end{align*}
by Lemma~\ref{lem:sn} and the initial assumption. Then this result and
the inequalities in (\ref{eq:tri-ineq-conj-1}) imply that posterior
consistency occurs if and only if
\[
E_M\left[\left.P_M\lp\left.\linf{\bntbg-\bnh}>\epsilon\given g,\sig,\bnh
,S_n\rp\given\bnh,S_n\right]\to0\,\aspo
\]
for every $\epsilon>0.$
Since $\bntbg-\bnh=(g+1)^{-1}(\gn-\bnh),$ we may equivalently state
that posterior consistency occurs if and only if $P_M[(g+1)^{-1}||\gn
-\bnh||_\infty>\epsilon\mid\bnh,S_n)\to0\,\aspo$ for every $\epsilon
>0.$ But again by the triangle inequality,
\begin{align}
&P_M\lp\left.\frac{1}{g+1}\linf{\gn-\bon}>2\epsilon\given\bnh,S_n\rp
-P_M\lp\left.\frac{1}{g+1}\linf{\bon-\bnh}>\epsilon\given\bnh,S_n\rp
\notag\\
&\le
P_M\lp\left.\frac{1}{g+1}\linf{\gn-\bnh}>\epsilon\given\bnh,S_n\rp
\label{eq:tri-ineq-conj-2}\\
&\le
P_M\lp\left.\frac{1}{g+1}\linf{\gn-\bon}>\epsilon/2\given\bnh,S_n\rp
\notag\\
&\qquad
+P_M\lp\left.\frac{1}{g+1}\linf{\bon-\bnh}>\epsilon/2\given\bnh,S_n\rp.
\notag
\end{align}
For any arbitrary $\epsilon>0,$
\begin{align*}
P_M\lp\left.\frac{1}{g+1}\linf{\bon-\bnh}>\epsilon\given\bnh,S_n\rp
&\le
P_M\lp\left.\linf{\bon-\bnh}>\epsilon\given\bnh,S_n\rp\\
&=
I\lp\linf{\bon-\bnh}>\epsilon\rp\to0\;\aspo
\end{align*}
by Lemma~\ref{lem:mle-linf}, where $I(\cdot)$ denotes the usual
indicator function. Then this result and~(\ref{eq:tri-ineq-conj-2})
together imply that
posterior consistency occurs if and only if $P_M[(g+1)^{-1}||\gn-\bon
||_\infty>\epsilon\mid\bnh,S_n]\to0\,\aspo$ for every $\epsilon>0.$
\end{prf}
%

%
\begin{prf}[Proof of Lemma~\ref{lem:pcg-hier-conj}]
From the form of the posterior in~(\ref{eq:post-g-conj}) and the
transformation in~(\ref{eq:transform}),
\begin{align*}
&\qquad\frac{T_n^2}{n^3}E_M\left[\left.\frac{g^2}{(g+1)^4}\given\bnh
,S_n\right]\\
&\le
\frac{T_n^2}{n^3}E_M\left[\left.\frac{1}{(g+1)^2}\given\bnh,S_n\right]
\\
&=
\frac{\displaystyle T_n^2\int_0^\infty(g+1)^{(n-p_n+a-c-6)/2}\left
[(g+1)(S_n+b)+T_n\right]^{-(n+a-2)/2}\;dg}{\displaystyle n^3\int
_0^\infty(g+1)^{(n-p_n+a-c-2)/2}\left[(g+1)(S_n+b)+T_n\right
]^{-(n+a-2)/2}\;dg} \\
&=
\frac{\displaystyle(S_n+b)^2\int_{W_n}^1 u^{(n-p_n+a-c-6)/2}
(1-u)^{(p_n+c)/2} \;du}
{\displaystyle n^3\int_{W_n}^1 u^{(n-p_n+a-c-2)/2} (1-u)^{(p_n+c-4)/2}
\;du}.
\end{align*}
Now let $H_n\sim\text{Beta}((n-p_n+a-c-4)/2,\;(p_n+c-2)/2)$ and
$\widetilde{H}_n\sim\text{Beta}((n-p_n+a-c)/2,\;(p_n+c-2)/2)$ with both
independent of $\bnh$ and $S_n$ under $P_M,$ and observe that $H_n$ is
stochastically smaller than $\widetilde{H}_n$ under~$P_M.$ Also let
$\Gamma(\cdot)$ denote the usual gamma function. Continuing, we have that
\begin{align*}
&\qquad\frac{T_n^2}{n^3}\;E_M\left[\left.\frac{g^2}{(g+1)^4}\given\bnh
,S_n\right]\\
&=
\frac{(S_n+b)^2\;\Gamma\lp\frac{n-p_n+a-c-4}{2}\rp\Gamma\lp\frac
{p_n+c+2}{2}\rp P_M\lp H_n>W_n\mid\bnh,S_n\rp}{n^3\;\Gamma\lp\frac
{n-p_n+a-c}{2}\rp\Gamma\lp\frac{p_n+c-2}{2}\rp P_M\lp\widetilde
{H}_n>W_n\mid\bnh,S_n\rp} \\
&\le
\frac{1}{n}\lp\frac{S_n+b}{n}\rp^2\frac
{(p_n+c)(p_n+c-2)}{(n-p_n+a-c-2)(n-p_n+a-c-4)}
\to0\;\;\aspo
\end{align*}
by Lemma~\ref{lem:sn}.
\end{prf}
%

%
\begin{prf}[Proof of Lemma~\ref{lem:wn-bound}]
Assume $\liminf_{n\to\infty}||\gn-\bon||_2^2\ge\delta$ for some $\delta>0$.
Then
\begin{align*}
&\qquad\limsup_{n\to\infty}W_n\\
&=\limsup_{n\to\infty}\lp1+\frac{T_n}{S_n+b}\rp^{-1} \\
&\le\lp1+\liminf_{n\to\infty}\left[\frac{p_n\sigo+n\lti\ltwo{\gn-\bon
}^2}{n-p_n}\right]\liminf_{n\to\infty}\left[\frac{(n-p_n)\;
T_n}{(S_n+b)\;\thon}\right]\rp^{-1} \\
&\le\lp1+\frac{\alpha}{1-\alpha}+\frac{\delta}{\lmax(1-\alpha)\sigo}\rp^{-1}
=\frac{(1-\alpha)\lmax\sigo}{\delta+\lmax\sigo}<1-\alpha
\;\;\aspo
\end{align*}
by Lemmas~\ref{lem:sn} and~\ref{lem:tn-to-mean}.
\end{prf}
%

%
\begin{prf}[Proof of Lemma~\ref{lem:wn-ln-to-zero}]
Assume $\ltwo{\gn-\bon}^2\to\infty,$ and let $\epsilon>0.$
Then
\begin{align*}
W_n
&=\lp1+\frac{T_n}{S_n+b}\rp^{-1}\\
&=\lp1+\frac{p_n\sigo+n\lti\ltwo{\gn-\bon}^2}{n-p_n}
\left[\frac{(n-p_n)\;T_n}{(S_n+b)\;\thon}\right]\rp^{-1} \to0\;\aspo
\end{align*}
since the term in square brackets converges to $1/\sigo\,\aspo$ by
Lemmas~\ref{lem:sn} and~\ref{lem:tn-to-mean}. This establishes~(i). Now write
\begin{align*}
\frac{\epsilon^{-1}\linf{\gn-\bon}(S_n+b)}{\epsilon^{-1}\linf{\gn-\bon
}(S_n+b)+T_n}
&=
\lp1+\frac{\epsilon\;T_n}{\linf{\gn-\bon}(S_n+b)}\rp^{-1}
\end{align*}
and observe that
\begin{align*}
&\qquad\lp1+\frac{\epsilon\;T_n}{\linf{\gn-\bon}(S_n+b)}\rp^{-1}\\
&=
\lp1+\frac{\epsilon\lp p_n\sigo+n\lti\ltwo{\gn-\bon}^2\rp}{(n-p_n)\linf
{\gn-\bon}}\left[\frac{(n-p_n)\;T_n}{(S_n+b)\;\thon}\right]\rp^{-1} \\
&\le
\lp1+\frac{\epsilon\ltwo{\gn-\bon}}{\lmax}\left[\frac{(n-p_n)\;
T_n}{(S_n+b)\;\thon}\right]\rp^{-1}\to0\;\aspo
\end{align*}
since, once again, the term in square brackets converges to $1/\sigo\,
\aspo$ by Lemmas~\ref{lem:sn} and~\ref{lem:tn-to-mean}. It then follows
immediately that $L_n(\epsilon)\to0\;\aspo,$ establishing~(ii).
\end{prf}
%

%
\begin{prf}[Proof of Lemma~\ref{lem:ln-bounds}]
Assume that $||\gn-\bon||_2^2\to A>0$ and $\liminf_{n\to\infty}||\gn
-\bon||_\infty>0.$ Let $\epsilon>0$.
Then $\limsup_{n\to\infty} W_n\le(1-\alpha)\lmax\sigo/(A+\lmax\sigo)<1\,
\aspo$ by Lemma~\ref{lem:wn-bound}, and
\begin{align*}
\limsup_{n\to\infty}\widetilde{L}_n(\epsilon)
&=
\limsup_{n\to\infty}\lp1+\frac{\epsilon\;T_n}{\linf{\gn-\bon}(S_n+b)}\rp
^{-1} \\
&=
\limsup_{n\to\infty}\lp1+\frac{\epsilon\lp p_n\sigo+n\lti\ltwo{\gn-\bon
}^2\rp}{(n-p_n)\linf{\gn-\bon}}\left[\frac{(n-p_n)\;T_n}{(S_n+b)\;\thon
}\right]\rp^{-1} \\
&\le\limsup_{n\to\infty}
\lp1+\frac{\epsilon\ltwo{\gn-\bon}}{\lmax}\left[\frac{(n-p_n)\;
T_n}{(S_n+b)\;\thon}\right]\rp^{-1} \\
&=\lp1+\frac{\epsilon\sqrt{A}}{\lmax\sigo}\rp^{-1}
=\frac{\lmax\sigo}{A^{1/2}\epsilon+\lmax\sigo}<1\;\;\aspo
\end{align*}
since the term in square brackets converges to $1/\sigo\,\aspo$ by
Lemmas~\ref{lem:sn} and~\ref{lem:tn-to-mean}. Define
\begin{align*}
L^\star(\epsilon)=\max\left\{\frac{(1-\alpha)\lmax\sigo}{A+\lmax\sigo
},\;\frac{\lmax\sigo}{A^{1/2}\epsilon+\lmax\sigo}\right\}<1,
\end{align*}
and observe that $\limsup_{n\to\infty}L_n(\epsilon)\le L^\star(\epsilon
)\,\aspo$.
This establishes~(i).

Now define $\widetilde{A}\defined\liminf_{n\to\infty}||\gn-\bon||_\infty
>0,$ and note that
\begin{align*}
&\qquad\liminf_{n\to\infty} L_n(\epsilon)\\
&\ge
\liminf_{n\to\infty} \lp1+\frac{\epsilon\;T_n}{\linf{\gn-\bon
}(S_n+b)}\rp^{-1} \\
&\ge
\lp1+\limsup_{n\to\infty}\left[\frac{\epsilon\lp p_n\sigo+n\lti\ltwo{\gn
-\bon}^2\rp}{(n-p_n)\linf{\gn-\bon}}\right]
\limsup_{n\to\infty}\left[\frac{(n-p_n)\;T_n}{(S_n+b)\;\thon}\right]\rp^{-1},
\end{align*}
which implies that
\begin{align*}
\liminf_{n\to\infty} L_n(\epsilon)
&\ge
\lp1+\frac{\epsilon\lp\alpha\sigo+A/\lmin\rp}{(1-\alpha)\widetilde
{A}\sigo}
\rp^{-1}\;\aspo
\end{align*}
by Lemmas~\ref{lem:sn} and~\ref{lem:tn-to-mean}.
Then it can be seen that for any $\zeta<1,$ there exists $\epsilon
_{\zeta}>0$ such that $\liminf_{n\to\infty} L_n(\epsilon_{\zeta})>\zeta
\,\aspo,$ establishing~(ii).
\end{prf}
%

%
\begin{prf}[Proof of Lemma~\ref{lem:beta-to-mean}]
Let $\epsilon>0$. Note that $E(Z_n)=a_n/(a_n+b_n)\to1-\alpha$, and thus
$|a_n/(a_n+b_n)-(1-\alpha)|\le\epsilon/2$ for all sufficiently large~$n$.
Also note that $\text{Var}(Z_n)=a_nb_n/[(a_n+b_n)^2(a_n+b_n+1)]\le
1/a_n<2/[n(1-\alpha)]$ for all sufficiently large~$n$. Then for all
sufficiently large~$n$,
\begin{align*}
&P\lp1-\alpha-\epsilon\le Z_n\le1-\alpha+\epsilon\rp\\
&\qquad=
P\lp1-\alpha-\frac{a_n}{a_n+b_n}-\epsilon\le Z_n-\frac{a_n}{a_n+b_n}\le
1-\alpha-\frac{a_n}{a_n+b_n}+\epsilon\rp\\
&\qquad\ge P\lp-\frac{\epsilon}{2}\le Z_n-\frac{a_n}{a_n+b_n}\le\frac
{\epsilon}{2}\rp
\ge
1-\frac{4}{\epsilon^2}\text{Var}(Z_n)
\ge
1-\frac{8}{n(1-\alpha)\epsilon^2}\to1,
\end{align*}
where the second of the three inequalities is Chebyshev's inequality.
\end{prf}
%

%
\begin{prf}[Proof of Theorem~\ref{thm:pcg-conj}]
By Lemmas~\ref{lem:pcg-hier} and~\ref{lem:pcg-hier-conj}, posterior
consistency occurs if and only if $P_M[(g+1)^{-1}||\gn-\bon||_\infty
>\epsilon\mid\bnh,S_n]\to0\,\aspo$ for every $\epsilon>0,$ which by
(\ref{eq:p-g-form}) occurs if and only if $P_M[W_n<U_n<L_n(\epsilon)\mid
\bnh,S_n]/P_M(U_n>W_n\mid\bnh,S_n)\to0\,\aspo$ for every $\epsilon>0.$
We now consider the same three cases as in the proof of Theorem~\ref
{thm:pcg-eb}.

Case~1: Suppose there do not exist a subsequence $k_n$ of $n$ and a
constant $A>0$ such that $||\gkn-\bokn||_2^2\to A$ and $||\gkn-\bokn
||_\infty\nrightarrow0.$ Let $k_n$ be a subsequence of $n$, and let
$\epsilon>0$. Now consider two sub-cases.

Case~1.1: Suppose $||\gkn-\bokn||_\infty\to0.$ Then $\epsilon^{-1}||\gkn
-\bokn||_\infty<1$ for all sufficiently large $n$. This implies that
$L_{k_n}(\epsilon)=W_{k_n}$ for all sufficiently large~$n\,\aspo$, and
therefore$P_M[W_{k_n}<U_{k_n}<L_{k_n}(\epsilon)\mid\bnh,S_n]=0$ for all
sufficiently large $n\,\aspo$. Also, $P_M(U_{k_n}>W_{k_n}\mid\bknh
,S_{k_n})>0$ for all~$n\,\aspo$ since $W_{k_n}<1$ for all~$n\,\aspo$. Thus,
\begin{align*}
\frac{P_M\left[W_{k_n}<U_{k_n}<L_{k_n}(\epsilon)\mid\bknh,S_{k_n}\right
]}{P_M(U_{k_n}>W_{k_n}\mid\bnh,S_n)}\to0\;\aspo
\end{align*}
by the combination of our results for its numerator and denominator.

Case~1.2: Note that Case~1.1 can be applied to any further
subsequence~$m_n$\break of $k_n$ for which $||\gmn-\bomn||_\infty\to0,$ so we
may suppose for Case~1.2 that\break $\liminf_{n\to\infty}||\gkn-\bokn||_\infty
>0.$ Note also that in this case, there cannot exist any further
subsequence~$m_n$ of~$k_n$ for which $||\gmn-\bomn||_2^2$ converges to
a nonzero constant, since this would contradict the original
supposition of Case~1. Then since $\liminf_{n\to\infty}||\gkn-\bokn
||_2^2\ge\liminf_{n\to\infty}||\gkn-\bokn||_\infty^2>0,$ it
follows\vadjust{\goodbreak}
that $||\gkn-\bokn||_2^2\to\infty.$ Then Lemma~\ref{lem:wn-ln-to-zero}
implies that both $W_{k_n}\to0\,\aspo$ and $L_n(\epsilon)\to0\,\aspo$,
which in turn implies that both $W_{k_n}<(1-\alpha)/2$ and
$L_{k_n}(\epsilon)<(1-\alpha)/2$ for all sufficiently large~$n\,\aspo$.
Then for all sufficiently large~$n\,\aspo$,
\begin{align*}
\frac{P_M\left[\left.W_{k_n}<U_{k_n}<L_{k_n}(\epsilon)\given\bknh
,S_{k_n}\right]}{P_M\lp U_{k_n}>W_{k_n}\mid\bknh,S_{k_n}\rp}
&\le
\frac{P_M\left[\left.U_{k_n}<\dfrac{1-\alpha}{2}\given\bknh
,S_{k_n}\right]}{P_M\left[\left.U_{k_n}>\dfrac{1-\alpha}{2}\given\bknh
,S_{k_n}\right]}
\to0\;\aspo
\end{align*}
by Lemma~\ref{lem:beta-to-mean}.
Finally, since Cases~1.1 and~1.2 together establish that the relevant
condition holds along any subsequence~$k_n$, it holds for the whole
sequence, and therefore posterior consistency occurs.

Case 2: Now suppose there exist a subsequence $k_n$ of $n$ and a
constant $A>0$ such that $||\gkn-\bokn||_2^2\to A>0$ and $||\gkn-\bokn
||_\infty\nrightarrow0,$ and suppose $\alpha=0.$
Note that Case~1 can be applied to any subsequence~$m_n$ of~$n$ for
which either $||\gmn-\bomn||_2^2$ does not converge to any nonzero
constant or $||\gmn-\bomn||_\infty\to0$, so it suffices to show that
the relevant condition holds along the subsequence~$k_n$. Note also
that this means we may suppose for Case~2 that $\liminf_{n\to\infty
}||\gkn-\bokn||_\infty>0$.
Now let \mbox{$\epsilon>0$}. By Lemma~\ref{lem:wn-bound}, $\limsup_{n\to\infty
}W_{k_n}\le\lmax\sigo/(A+\lmax\sigo)\,\aspo$, which implies that
$W_{k_n}<2\lmax\sigo/(A+2\lmax\sigo)$ for all sufficiently large~$n\,
\aspo$. Moreover, by Lemma~\ref{lem:ln-bounds}, there exists $L^\star
(\epsilon)<1$ such that $\limsup_{n\to\infty}L_{k_n}(\epsilon)\le
L^\star(\epsilon)\,\aspo$, which implies that $L_{k_n}(\epsilon
)<[1+L^\star(\epsilon)]/2$ for all sufficiently large~$n\,\aspo$. Then
for all sufficiently large~$n\,\aspo$,
\begin{align*}
&\qquad\frac{P_M\left[\left.W_{k_n}<U_{k_n}<L_{k_n}(\epsilon)\given
\bknh,S_{k_n}\right]}{P_M\lp\left. U_{k_n}>W_{k_n}\given\bknh,S_{k_n}\rp
}\\
&\le
\frac{P_M\left[\left.U_{k_n}<\dfrac{1+L^\star(\epsilon)}{2}\given\bknh
,S_{k_n}\right]}{P_M\lp\left.U_{k_n}>\dfrac{2\lmax\sigo}{A+2\lmax\sigo
}\given\bknh,S_{k_n}\rp}\to0\;\aspo
\end{align*}
by Lemma~\ref{lem:beta-to-mean}. Therefore posterior consistency occurs.

Case 3: Now suppose there exist a subsequence $k_n$ of $n$ and a
constant $A>0$ such that $||\gkn-\bokn||_2^2\to A>0$ and $||\gkn-\bokn
||_\infty\nrightarrow0,$ but suppose $\alpha>0.$
By Lemma~\ref{lem:wn-bound}, $\limsup_{n\to\infty}W_{k_n}\le(1-\alpha
)\lmax\sigo/(A+\lmax\sigo)\,\aspo$, which implies that
$W_{k_n}<2(1-\alpha)\lmax\sigo/(A+2\lmax\sigo)$ for all sufficiently
large~$n\,\aspo$.
By Lemma~\ref{lem:ln-bounds}, there exists $\epsilon_{1-\alpha/4}>0$
such that
$\liminf_{n\to\infty}L_{k_n}(\epsilon_{1-\alpha/4})\ge1-\alpha/4\,\aspo
$, which implies that $L_{k_n}(\epsilon_{1-\alpha/4})>1-\alpha/2$ for
all sufficiently large~$n\,\aspo$. Then for all sufficiently large~$n\,
\aspo$,
\begin{align*}
&\qquad\frac{P_M\left[\left.W_{k_n}<U_{k_n}<L_{k_n}(\epsilon_{1-\alpha
/4})\given\bknh,S_{k_n}\right]}{P_M\lp\left. U_{k_n}>W_{k_n}\given\bknh
,S_{k_n}\rp}\\
&\ge
P_M\left[\left.W_{k_n}<U_{k_n}<L_{k_n}(\epsilon_{1-\alpha/4})\given\bknh
,S_{k_n}\right]\\
&\ge
P_M\left[\left.\frac{2(1-\alpha)\lmax\sigo}{A+2\lmax\sigo
}<U_{k_n}<1-\frac{\alpha}{2}\given\bknh,S_{k_n}\right]\\
&\to1\;\aspo
\end{align*}
by Lemma~\ref{lem:beta-to-mean}. Since the relevant condition fails to
hold for the subsequence~$k_n$, it fails to hold for the overall
sequence. Therefore posterior consistency does not occur.
\end{prf}

%
\begin{prf}[Proof of Lemma~\ref{lem:pcg-hier-zs}]
Consider two cases.

Case~1: Suppose $||\gn-\bon||_2^2\to0.$ Then $\thon=p_n\sigo+n\lton||\gn
-\bon||_2^2\le n\sigo$ for all sufficiently large~$n.$ This result and
(\ref{eq:nccs-fourth}) imply that
%
\begin{align*}
\lp\mu_4\rp_0\lp T_n\rp
&\defined
E_0\left[\lp T_n-\thon\rp^4\right]
\le
48\sigma_0^4\thon^2+192\sigma_0^6\thon\le96n^2\sigma_0^8
\end{align*}
for all sufficiently large~$n$.
Then there exists $N$ such that
\begin{align*}
\sum_{n=N}^\infty P_0\lp T_n>2n\sigo\rp
&\le
\sum_{n=N}^\infty P_0\lp\left|T_n-\thon\right|>n\sigo\rp
\le
\sum_{n=N}^\infty\frac{96n^2\sigma_0^8}{n^4\sigma_0^8}=96\sum
_{n=N}^\infty\frac{1}{n^2}<\infty
\end{align*}
by Markov's inequality applied to $(T_n-\thon)^4,$ which in turn
implies by the Borel-Cantelli lemma that $\limsup_{n\to\infty}
(T_n/n)\le2\sigo\,\aspo.$ Therefore,
$n^{-3}\,T_n^2\,E_M[g^2(g+1)^{-4}\mid\bnh,S_n]\le n^{-3}\,T_n^2\to0\,
\aspo.$

Case~2: Note immediately that Case~1 can be applied to any
subsequence~$k_n$ of~$n$ for which $||\gn-\bon||_2^2\to0,$ so we may
suppose for Case~2 that $\liminf_{n\to\infty}||\gn-\bon||_2^2>0.$ Then
$\limsup_{n\to\infty} W_n<1-\alpha$ by Lemma~\ref{lem:wn-bound}. Define
$
\psi_n(u)\defined I_{(W_n,1)}(u)\exp[-nW_n(1-u)/2(u-W_n)],
$
where $I$ denotes the usual indicator function,
and note that this is a nondecreasing function of $u$ on the interval~$(0,1).$
Using the form of the posterior in~(\ref{eq:post-g-zs}) and the
transformation in~(\ref{eq:transform}), we may write
\begin{align*}
&\frac{T_n^2}{n^3}E_M\left[\left.\frac{g^2}{(g+1)^4}\given\bnh,S_n\right
] \\
&\qquad=
\frac{\displaystyle T_n^2\int_0^\infty\frac{(g+1)^{(n-p_n+a-10)/2}}{\big
[(g+1)(S_n+b)+T_n\big]^{-(n+a-2)/2}}\;
g^{1/2}\exp\lp-\frac{n}{2g}\rp\;dg}{\displaystyle n^3\int_0^\infty\frac
{(g+1)^{(n-p_n+a-2)/2}}{\big[(g+1)(S_n+b)+T_n\big]^{-(n+a-2)/2}}\;
g^{-3/2}\exp\lp-\frac{n}{2g}\rp\;dg} \\
&\qquad=
\frac{\displaystyle(S_n+b)^4\int_0^1 u^{(n-p_n+a-c-10)/2}
(1-u)^{(p_n+4)/2}\left[\frac{u-W_n}{W_n(1-u)}\right]^{1/2}\psi_n(u)\;du}
{\displaystyle n^3\; T_n^2\int_0^1 u^{(n-p_n+a-c-2)/2}
(1-u)^{(p_n-4)/2}\left[\frac{u-W_n}{W_n(1-u)}\right]^{-3/2}\psi_n(u) \;
du} \\
&\qquad\le
\frac{\displaystyle(S_n+b)^2\int_0^1 u^{(n-p_n+a-c-9)/2}
(1-u)^{(p_n+3)/2}\;\psi_n(u) \;du}
{\displaystyle n^3 (1-W_n)^2\int_0^1 u^{(n-p_n+a-c-5)/2}
(1-u)^{(p_n-1)/2}\;\psi_n(u) \;du}.
\end{align*}
Now let $h_n$ and $\widetilde{h}_n$ denote the densities with respect
to Lebesgue measure of $\text{Beta}((n-p_n+a-7)/2,\;(p_n+5)/2)$ and
$\text{Beta}((n-p_n+a-3)/2,\;(p_n+1)/2)$ random variables,
respectively. Then we may continue by writing
\begin{align*}
&\qquad\frac{T_n^2}{n^3}E_M\left[\left.\frac{g^2}{(g+1)^4}\given\bnh
,S_n\right]\\
&\le
\frac{\displaystyle(S_n+b)^2\;\Gamma\lp\frac{n-p_n+a-7}{2}\rp\Gamma\lp
\frac{p_n+5}{2}\rp\;\int_0^1 h_n(u)\;\psi_n(u) \;du}
{\displaystyle n^3 (1-W_n)^2\;\Gamma\lp\frac{n-p_n+a-3}{2}\rp\Gamma\lp
\frac{p_n+1}{2}\rp\;\int_0^1 \widetilde{h}_n(u)\;\psi_n(u) \;du} \\
&\le
\frac{(S_n+b)^2\;(p_n+3)(p_n+1)}
{n^3 (1-W_n)^2\;(n-p_n+a-5)(n-p_n+a-7)}\to0\;\aspo.
\end{align*}
Note that the last inequality holds because a random variable with
density $h_n$ is stochastically smaller than a random variable with
density $\widetilde{h}_n$ and because $\psi_n$ is nondecreasing on
$(0,1),$ while the almost sure convergence to zero is by Lemma~\ref
{lem:sn} and the fact that $\limsup_{n\to\infty} W_n<1-\alpha\le1\,\aspo
$ by Lemma~\ref{lem:wn-bound}.
\end{prf}
%

%
\begin{prf}[Proof of Lemma~\ref{lem:beta-tail}]
Note immediately that both (i) and (ii) are trivial if $\xi=0$ or $\xi
\ge1$, so assume $0<\xi<1$.
Next, by Stirling's approximation, we may bound the normalizing
constant by
\begin{align*}
\log\frac{\Gamma(a_n+b_n)}{\Gamma(a_n)\Gamma(b_n)}
&\le
\log\frac{(a_n+b_n)^{a_n+b_n-1/2}}{a_n^{a_n-1/2}\;b_n^{b_n-1/2}}
\end{align*}
for all sufficiently large~$n$. We may rewrite this as
\begin{align*}
\log\frac{\Gamma(a_n+b_n)}{\Gamma(a_n)\Gamma(b_n)}
&\le a_n\log\lp\frac{a_n+b_n}{a_n}\rp+b_n\log\lp\frac{a_n+b_n}{b_n}\rp
+\frac{1}{2}\log\lp\frac{a_n b_n}{a_n+b_n}\rp
\end{align*}
for all sufficiently large~$n$. Then
\begin{align*}
&\qquad P(Z_n\le\xi)\\
&=
\frac{\Gamma(a_n+b_n)}{\Gamma(a_n)\Gamma(b_n)}\int_0^\xi
z^{a_n-1}(1-z)^{b_n-1}\;dz
\\
&\le
\frac{\Gamma(a_n+b_n)}{\Gamma(a_n)\Gamma(b_n)}\int_0^\xi z^{a_n-1}\;dz
=
\frac{\Gamma(a_n+b_n)\;\xi^{a_n}}{\Gamma(a_n)\Gamma(b_n)\;a_n} \\
&\le
\exp\left[a_n\log\xi -\log a_n + a_n\log\lp\frac{a_n+b_n}{a_n}\rp
+b_n\log\lp\frac{a_n+b_n}{b_n}\rp+\frac{1}{2}\log\lp\frac{a_n
b_n}{a_n+b_n}\rp\right]
\end{align*}
for all sufficiently large~$n.$ Now observe that
\begin{align*}
&\qquad\frac{1}{n}\log P(Z_n\le\xi)\\
&\le
\frac{a_n}{n}\log\xi -\frac{1}{n}\log a_n + \frac{a_n}{n}\log\lp\frac
{a_n+b_n}{a_n}\rp
+\frac{b_n}{n}\log\lp\frac{a_n+b_n}{b_n}\rp+\frac{1}{2n}\log\lp\frac
{a_n b_n}{a_n+b_n}\rp\\
&\to
\begin{cases}
(1-\alpha)\log\xi-(1-\alpha)\log(1-\alpha)-\alpha\log\alpha& \text
{if } \alpha>0,\\
\log\xi &\text{if } \alpha=0.
\end{cases}
\end{align*}
If $\alpha>0$, then $(1-\alpha)\log(1-\alpha)+\alpha\log\alpha\ge-\log
2$, and thus $\limsup_{n\to\infty} n^{-1}\log P(Z_n\le\xi)\le(1-\alpha
)\log\xi+\log2$. Then $n^{-1}\log P(Z_n\le\xi)\le(1-\alpha)\log\xi
+\log4$ for all sufficiently large~$n,$ which implies~(i). If instead
$\alpha=0$, then $\limsup_{n\to\infty} n^{-1}\log P(Z_n\le\xi)\le\log
\xi,$ so $n^{-1}\log P(Z_n\le\xi)\le\frac{1}{2}\log\xi$ for all
sufficiently large~$n$ (noting that $\log\xi<0$). This implies~(ii).
\end{prf}
%

%
\begin{prf}[Proof of Lemma~\ref{lem:rn-bound}]
Let $\delta=\liminf_{n\to\infty}||\gn-\bon||_2^2>0$. Then by Lemma~\ref
{lem:wn-bound}, $\limsup_{n\to\infty}W_n\le(1-\alpha)\lmax\sigo/(\delta
+\lmax\sigo)\,\aspo$, which implies that $W_n<2(1-\alpha)\lmax\sigo
/(\delta+2\lmax\sigo)<1-\alpha$ for all sufficiently large~$n\,\aspo$.
Then for all sufficiently large~$n\,\aspo$,
\begin{align}
R_n
&\ge
\int_{(1-\alpha+W_n)/2}^1 f_n(u)\left[\frac{u-W_n}{(1-u)}\right]^{-3/2}
\exp\left[-\frac{nW_n(1-u)}{2(u-W_n)}\right]\;du \notag\\[3pt]
&\ge
\exp\lp-\frac{nW_n}{1-\alpha-W_n}\rp
\int_{(1-\alpha+W_n)/2}^1 f_n(u)\left[\frac{u}{(1-u)}\right
]^{-3/2}du\notag\\[3pt]
&=
\frac{\Gamma\lp\dfrac{n-p_n+a-3}{2}\rp\Gamma\lp\dfrac{p_n+1}{2}\rp}
{\Gamma\lp\dfrac{n-p_n+a}{2}\rp\Gamma\lp\dfrac{p_n-2}{2}\rp}
\;\exp\lp-\frac{nW_n}{1-\alpha-W_n}\rp\notag\\[3pt]
&\qquad\times P_M\lp\left.\frac{1-\alpha+W_n}{2}<\widetilde
{U}_n<1\given\bnh,S_n\rp\notag\\[3pt]
&\ge
\frac{\Gamma\lp\dfrac{n-p_n+a-3}{2}\rp\Gamma\lp\dfrac{p_n+1}{2}\rp}
{\Gamma\lp\dfrac{n-p_n+a}{2}\rp\Gamma\lp\dfrac{p_n-2}{2}\rp}
\;\exp\lp-\frac{nW_n}{1-\alpha-W_n}\rp\notag\\[3pt]
&\qquad\times
P_M\left[\lp\frac{\delta+4\lmax\sigo}{2\delta+4\lmax\sigo}\rp\lp1-\alpha
\rp<\widetilde{U}_n<1\right]
\label{eq:rn-lower}
\end{align}
where $\widetilde{U}_n\sim\text{Beta}((n-p_n+a-3)/2,(p_n+1)/2),$
independent of $\bnh$ and $S_n$, under~$P_M$.
For all sufficiently large~$n,$ Stirling's approximation yields
that\vadjust{\eject}
\begin{align*}
&\frac{\Gamma\lp\dfrac{n-p_n+a-3}{2}\rp\Gamma\lp\dfrac{p_n+1}{2}\rp}
{\Gamma\lp\dfrac{n-p_n+a}{2}\rp\Gamma\lp\dfrac{p_n-2}{2}\rp}\\
&\ge
\frac{
\lp\dfrac{n-p_n+a-3}{2}\rp^{(n-p_n+a-4)/2}\exp\lp-\dfrac
{n-p_n+a-3}{2}\rp}
{2\lp\dfrac{n-p_n+a}{2}\rp^{(n-p_n+a-1)/2}\exp\lp-\dfrac{n-p_n+a}{2}\rp
}\\
&\qquad\times\frac{\lp\dfrac{p_n+1}{2}\rp^{p_n/2}\exp\lp-\dfrac
{p_n+1}{2}\rp}
{\lp\dfrac{p_n-2}{2}\rp^{(p_n-3)/2}\exp\lp-\dfrac{p_n-2}{2}\rp}\\
&=
\frac{1}{2}
\lp\frac{n-p_n+a-3}{n-p_n+a}\rp^{(n-p_n+a-4)/2}\!
\lp\frac{p_n+1}{p_n-2}\rp^{(p_n-3)/2}\!
\lp\frac{p_n+1}{n-p_n+a}\rp^{3/2}\!.
\end{align*}
Then for all sufficiently large~$n$,
\begin{align}
\frac{\Gamma\lp\dfrac{n-p_n+a-3}{2}\rp\Gamma\lp\dfrac{p_n+1}{2}\rp}
{\Gamma\lp\dfrac{n-p_n+a}{2}\rp\Gamma\lp\dfrac{p_n-2}{2}\rp}
&\ge
\frac{1}{4}\lp\frac{p_n+1}{n-p_n+a}\rp^{3/2}
\ge2\,(4n)^{-3/2}.
\label{eq:rn-stirling}
\end{align}
Now observe that
\begin{align*}
P_M\left[\lp\frac{\delta+4\lmax\sigo}{2\delta+4\lmax\sigo}\rp\lp1-\alpha
\rp<\widetilde{U}_n<1\right]\to1
\end{align*}
by Lemma~\ref{lem:beta-to-mean}, which implies that
\begin{align}
P_M\left[\lp\frac{\delta+4\lmax\sigo}{2\delta+4\lmax\sigo}\rp\lp1-\alpha
\rp<\widetilde{U}_n<1\right]>\frac{1}{2}
\label{eq:rn-beta}
\end{align}
for all sufficiently large~$n$. Then by combining Inequalities~\ref
{eq:rn-lower},~\ref{eq:rn-stirling},~and~\ref{eq:rn-beta}, we have that
for all sufficiently large~$n\,\aspo$,
\begin{align*}
R_n
&\ge(4n)^{-3/2}\exp\lp-\frac{nW_n}{1-\alpha-W_n}\rp
=\exp\left\{-n\left[\frac{W_n}{1-\alpha-W_n}+\frac{3}{2n}\log(4n)\right
]\right\}.
\end{align*}
Finally, take $K=2\limsup_{n\to\infty}[W_n/(1-\alpha-W_n)]$. Observe
that $K<\infty\,\aspo$ due to the fact that $\limsup_{n\to\infty} W_n\le
(1-\alpha)\lmax\sigo/(\delta+\lmax\sigo)<1-\alpha\,\aspo.$
Then $R_n\ge\exp(-nK)$ for all sufficiently large~$n\,\aspo$.
\end{prf}
%

%
\begin{prf}[Proof of Lemma~\ref{lem:qn-bound}]
Let $\epsilon>0,$ and assume $||\gn-\bon||_2^2\to\infty.$
Then by Lemma~\ref{lem:wn-ln-to-zero}, $W_n\to0\,\aspo$ and
$L_n(\epsilon)\to0\,\aspo.$
Next, observe that the last two terms of the integrand in $Q_n(\epsilon
)$ comprise an unnormalized $\text{InverseGamma}(1/2,\;nW_n/2)$ density
in $(u-W_n)/(1-u),$ the mode of which occurs at $nW_n/3.$
Then for all sufficiently large~$n,$
\begin{align*}
Q_n(\epsilon)
\le
\int_{W_n}^{L_n(\epsilon)}f_n(u)
\lp\frac{nW_n}{3}\rp^{-3/2}\exp\lp-\frac{3}{2}\rp\;du
&\le
2\lp nW_n\rp^{-3/2}\int_{0}^{L_n(\epsilon)} f_n(u)\;du \\
&\le
2^{2n+1}\lp nW_n\rp^{-3/2}\left[L_n(\epsilon)\right]^{n(1-\alpha)}
\end{align*}
by Lemma~\ref{lem:beta-tail}.
Now note that if $\epsilon^{-1}||\gn-\bon||_\infty\le1,$ then
$L_n(\epsilon)=W_n,$ in which case $Q_n(\epsilon)=0$ and the result is
trivial. So instead assume that $\epsilon^{-1}||\gn-\bon||_\infty>1,$
which in turn implies that $L_n(\epsilon)\le\epsilon^{-1}||\gn-\bon
||_\infty W_n.$ Then
\begin{align*}
Q_n(\epsilon)
&\le
2^{2n+1}\lp\frac{n(S_n+b)}{S_n+b+T_n}\rp^{-3/2}\lp\frac{\epsilon
^{-1}||\gn-\bon||_\infty(S_n+b)}{\epsilon^{-1}||\gn-\bon||_\infty
(S_n+b)+T_n}\rp^{n(1-\alpha)}\\
&\le
2^{2n+1}\;n^{-3/2}
\lp1+\frac{\lp p_n\sigo+n\lti\ltwo{\gn-\bon}^2\rp}{n-p_n}\left[\frac
{(n-p_n)\;T_n}{(S_n+b)\;\thon}\right]\rp^{3/2} \\
&\qquad\times\lp1+\frac{\epsilon\lp p_n\sigo+n\lti\ltwo{\gn-\bon}^2\rp
}{(n-p_n)\linf{\gn-\bon}}\left[\frac{(n-p_n)\;T_n}{(S_n+b)\;\thon}\right
]\rp^{-n(1-\alpha)} \\
&\le
2^{2n+1}\;n^{-3/2}
\lp\frac{4\ltwo{\gn-\bon}^2}{(1-\alpha)\lton\sigo}\rp^{3/2}
\lp\frac{\epsilon\ltwo{\gn-\bon}}{2(1-\alpha)\lton\sigo}\rp^{-n(1-\alpha
)} \\
&=
2^{n(3-\alpha)+4}(n\epsilon)^{-3/2}\lp\epsilon^{-1}(1-\alpha)\lmax\sigo
\rp^{n(1-\alpha)-3/2}\ltwo{\gn-\bon}^{-n(1-\alpha)+3}
\end{align*}
for all sufficiently large~$n\,\aspo$ by Lemmas~\ref{lem:sn} and~\ref
{lem:tn-to-mean} since the quantity in square brackets converges to
$1/\sigo\,\aspo.$ Now continue by writing that for all sufficiently
large~$n\,\aspo,$
\begin{align*}
&\qquad Q_n(\epsilon)\\
&\le
\exp\left\{-n\left[
\lp1-\alpha-\frac{3}{n}\rp\log\lp\ltwo{\gn-\bon}\rp-\lp1-\alpha-\frac
{3}{2n}\rp\log\lp\epsilon^{-1}(1-\alpha)\lmax\sigo\rp\right.\right.\\
&\qquad\qquad\qquad\left.\left.{}+\frac{3}{2n}\log(n\epsilon)
-\lp3-\alpha+\frac{4}{n}\rp\log2
\right]\right\} \\
&=
\exp\left[-n\kappa_n(\epsilon)\right],
\end{align*}
where $\kappa_n(\epsilon)\to\infty$ is defined to be the quantity in
square brackets.
\end{prf}
%

%
\begin{prf}[Proof of Lemma~\ref{lem:zs-alpha-zero}]
Assume $||\gn-\bon||_2^2\to A>0,$ $\liminf_{n\to\infty}||\gn-\bon
||_\infty>0,$ and $\alpha=0.$ Let $\epsilon>0$. Note that $R_n>0$ for
all~$n\,\aspo$ since $W_n<1$ for all~$n\,\aspo$. Then whenever
$L_n(\epsilon)\le W_n$, we immediately have that $Q_n(\epsilon)/R_n=0$
exactly, so we may instead assume that $L_n(\epsilon)>W_n$ for all~$n$.
By Lemma~\ref{lem:ln-bounds}, there exists $L^\star(\epsilon)<1$ such
that $\limsup_{n\to\infty}L_n(\epsilon)\le
L^\star(\epsilon)\,\aspo$, which implies that $L_n(\epsilon)<[1+L^\star
(\epsilon)]/2$ for all sufficiently large~$n\,\aspo$.
Then we may write that for all sufficiently large~$n\,\aspo$,
\begin{align}
R_n
&\ge
\int_{[1+L_n(\epsilon)]/2}^1 f_n(u)\left[\frac{u-W_n}{(1-u)}\right]^{-3/2}
\exp\left[-\frac{nW_n(1-u)}{2(u-W_n)}\right]\;du \notag\\
&\ge
\exp\left\{-\frac{nW_n\left[1-L_n(\epsilon)\right]}{2\left
[1+L_n(\epsilon)-2W_n\right]}\right\}
\int_{[3+L^\star(\epsilon)]/4}^1 f_n(u)\left[\frac{u}{(1-u)}\right
]^{-3/2}du\notag\\
&\ge
\frac{\Gamma\lp\dfrac{n-p_n+a-3}{2}\rp\Gamma\lp\dfrac{p_n+1}{2}\rp}
{\Gamma\lp\dfrac{n-p_n+a}{2}\rp\Gamma\lp\dfrac{p_n-2}{2}\rp}
\;\exp\left\{-\frac{nW_n\left[1-L_n(\epsilon)\right]}{4\left
[L_n(\epsilon)-W_n\right]}\right\} \notag\\
&\qquad\times
P_M\lp\left.\frac{3+L^\star(\epsilon)}{4}<\widetilde{U}_n<1\given\bnh
,S_n\rp\notag\\
&\ge
(4n)^{-3/2}
\;\exp\left\{-\frac{nW_n\left[1-L_n(\epsilon)\right]}{4\left
[L_n(\epsilon)-W_n\right]}\right\}
\label{eq:zero-rn}
\end{align}
by Inequalities~\ref{eq:rn-stirling}~and~\ref{eq:rn-beta}.
Next, write $Q_n(\epsilon)$ as
\begin{align*}
Q_n(\epsilon)
&=
\int_{W_n}^{L_n(\epsilon)} f_n(u)\left[\frac{u-W_n}{(1-u)}\right]^{-3/2}
\exp\left[-\frac{W_n(1-u)}{2(u-W_n)}\right]
\exp\left[-\frac{(n-1)W_n(1-u)}{2(u-W_n)}\right]\;du.
\end{align*}
The second and third terms of the integrand comprise an unnormalized
\[
\text{InverseGamma}(1/2,W_n/2)
\]

\noindent
density in $(u-W_n)/(1-u)$, which has mode $W_n/3.$ Then
\begin{align*}
Q_n(\epsilon)
&\le
\lp\frac{W_n}{3}\rp^{-3/2}\;\exp\lp-\frac{3}{2}\rp\;
\exp\left\{-\frac{(n-1)W_n\left[1-L_n(\epsilon)\right]}{2\left
[L_n(\epsilon)-W_n\right]}\right\}
\int_{0}^{L_n(\epsilon)} f_n(u)\;du, \\
&\le
\lp2 W_n\rp^{-3/2}\exp\left\{-\frac{(n-1)W_n\left[1-L_n(\epsilon)\right
]}{2\left[L_n(\epsilon)-W_n\right]}\right\}
\left[L_n(\epsilon)\right]^{n/2}
\end{align*}
by Lemma~\ref{lem:beta-tail}.
Then this result and Inequality~\ref{eq:zero-rn} together yield that
for all sufficiently large~$n\,\aspo,$
\begin{align}
\frac{Q_n(\epsilon)}{R_n}
&\le
\lp\frac{W_n}{2n}\rp^{-3/2}\left[L_n(\epsilon)\right]^{n/2}
\exp\left\{
-\frac{nW_n\left[1-L_n(\epsilon)\right]}{2\left[L_n(\epsilon)-W_n\right
]}\lp\frac{n-1}{n}-\frac{1}{2}\rp\right\}\notag\\
&\le
\lp\frac{2n}{W_n}\rp^{3/2}
\exp\left\{-\frac{nW_n\left[1-L^\star(\epsilon)\right]}{16}\right\}.
\label{eq:zero-qr}
\end{align}
Now observe that
\begin{align*}
&\qquad\liminf_{n\to\infty}W_n\\
&= \liminf_{n\to\infty}\lp1+\frac{T_n}{S_n+b}\rp^{-1}\\
&=\liminf_{n\to\infty}\lp1+\frac{p_n\sigo+n\lti\ltwo{\gn-\bon}^2}{n-p_n}
\left[\frac{(n-p_n)\;T_n}{(S_n+b)\;\thon}\right]\rp^{-1} \\
&\ge\lp1+A/\lmin\sigo\rp^{-1}=\frac{\lmin\sigo}{A+\lmin\sigo}\;\;\aspo,
\end{align*}
which implies that $W_n>\lmin\sigo/(2A+\lmin\sigo)$ for all
sufficiently large~$n\,\aspo$. We may combine this with Inequality~\ref
{eq:zero-qr} to yield that for all sufficiently large~$n\,\aspo$,
\begin{align*}
\frac{Q_n(\epsilon)}{R_n}
&\le
\left[\frac{2n(2A+\lmin\sigo)}{\lmin\sigo}\right]^{3/2}
\exp\left\{-\frac{n\lmin\sigo\left[1-L^\star(\epsilon)\right
]}{16(2A+\lmin\sigo)}\right\}\to0\;\aspo
\end{align*}
since $L^\star(\epsilon)<1$.
\end{prf}

%
\begin{prf}[Proof of Theorem~\ref{thm:pcg-zs}]
By Lemmas~\ref{lem:pcg-hier} and~\ref{lem:pcg-hier-zs}, posterior
consistency occurs if $P_M[(g+1)^{-1}||\gn-\bon||_\infty>\epsilon\mid
\bnh,S_n]\to0\,\aspo$ for every $\epsilon>0.$ Then by~(\ref
{eq:g-zs-ratio}), this occurs if $Q_n(\epsilon)/R_n\to0\,\aspo$ for
every $\epsilon>0.$
We now proceed according to cases similar to those in the proofs of the
previous theorems.

Case~1: Suppose there do not exist a subsequence $k_n$ of $n$ and a
constant $A>0$ such that $||\gkn-\bokn||_2^2\to A$ and $||\gkn-\bokn
||_\infty\nrightarrow0.$ Let $k_n$ be a subsequence of $n$, and let
$\epsilon>0$. Now consider two sub-cases.

Case~1.1: Suppose $||\gkn-\bon||_\infty\to0.$ Then $\epsilon^{-1}||\gkn
-\bokn||_\infty<1$ for all sufficiently large~$n\,\aspo$. This implies
that $L_{k_n}(\epsilon)=W_{k_n}$ and $Q_{k_n}(\epsilon)=0$ for all
sufficiently large~$n\,\aspo$. Also, $R_{k_n}>0$ for all $n\,\aspo$
since $W_{k_n}<1\,\aspo$. Therefore, $Q_{k_n}(\epsilon)/R_{k_n}\to0\,
\aspo$.

Case~1.2: Note that Case~1.1 can be applied to any further
subsequence~$m_n$ of~$k_n$ for which $||\gmn-\bomn||_\infty\to0,$ so we
may suppose for Case~1.2 that $\liminf_{n\to\infty}||\gkn-\bokn||_\infty
>0.$ Note also that in this case, there cannot exist any further
subsequence~$m_n$ of~$k_n$ for which $||\gmn-\bomn||_2^2$ converges to
a nonzero constant, since this would contradict the original
supposition of Case~1. Then since $\liminf_{n\to\infty}||\gkn-\bokn
||_2^2\ge\liminf_{n\to\infty}||\gkn-\bokn||_\infty^2>0,$ it follows
that $||\gkn-\bokn||_2^2\to\infty.$ Observe that by Lemmas~\ref
{lem:rn-bound} and~\ref{lem:qn-bound}, there exist a constant~$K$ and a
sequence of constants $\kappa_n(\epsilon)\to\infty$ such that
$Q_{k_n}(\epsilon)/R_{k_n}\le\exp\left\{-n\left[\kappa_n(\epsilon
)-K\right]\right\}\to0\,\aspo$.
Finally, since Cases~1.1~and~1.2 together establish that
$Q_{k_n}(\epsilon)/R_{k_n}\to0\,\aspo$ for every subsequence~$k_n$, it
follows that $Q_n(\epsilon)/R_n\to0\,\aspo$, and therefore posterior
consistency occurs.

Case 2: Now suppose there exist a subsequence $k_n$ of $n$ and a
constant $A>0$ such that $||\gkn-\bokn||_2^2\to A>0$ and $||\gkn-\bokn
||_\infty\nrightarrow0,$ and suppose $\alpha=0.$
Note that Case~1 can be applied to any subsequence~$m_n$ of~$n$ for
which either $||\gmn-\bomn||_2^2$ does not converge to any nonzero
constant or $||\gmn-\bomn||_\infty\to0$, so it suffices to show that
$Q_{k_n}(\epsilon)/R_{k_n}\to0\,\aspo$. Note also that this means we
may suppose for Case~2 that $\liminf_{n\to\infty}||\gkn-\bokn||_\infty
>0$. Now let $\epsilon>0$.
Then we immediately have that $Q_{k_n}(\epsilon)/R_n\to0\,\aspo$ by
Lemma~\ref{lem:zs-alpha-zero}.
Therefore posterior consistency occurs.
\end{prf}

\bibliographystyle{ba}

\begin{acknowledgement}
This work was supported in part by NSF Grants DMS--1106084 and DMS--1007494.
\end{acknowledgement}

\end{document}